\documentclass[onecolumn, 12pt]{article}

\usepackage[utf8]{inputenc} 
\usepackage[T1]{fontenc}    
\usepackage{hyperref}       
\usepackage{url,mathtools}            
\usepackage{booktabs}       
\usepackage{amsfonts}       
\usepackage{nicefrac}       
\usepackage{microtype}      
\usepackage{graphicx}
\usepackage{color}
\usepackage{bm, bbm}

\usepackage{algorithm,algorithmic}
\usepackage{amsthm,amsmath,amssymb}
\usepackage{natbib}
\newcounter{thm}
\newcounter{re}
 \def\E{{\rm E}}
\def\Var{{\rm Var}}
\def\Cov{{\rm Cov}}

\def\tr{{\rm tr}}

 \newtheorem{theorem}[thm]{Theorem}
 \newtheorem{lemma}[thm]{Lemma}
 
 \newtheorem{proposition}[thm]{Proposition}

\title{On near-optimality of one-sample update for sequential joint detection and estimation}

\author{
 Yang Cao, Yao Xie, and Huan Xu
 \\
  H. Milton Stewart School of Industrial and Systems Engineering\\
  Georgia Institute of Technology\\
  Atlanta, GA 30342 \\
  \texttt{\{caoyang, yao.xie, huan.xu\}@isye.gatech.edu} \\
}

\begin{document}

\title{$\textsf{S}^3 \textsf{T}$: An Efficient Score Statistic for Spatio-Temporal Surveillance}

\author{Junzhuo Chen, Seong-Hee Kim and Yao Xie\\H. Milton Stewart School of Industrial Engineering\\ Georgia Institute of Technology}

\maketitle

\begin{abstract}
	
We present an efficient score statistic, called the $\textsf{S}^3 \textsf{T}$ statistic, to detect the emergence of a spatially and temporally correlated signal from either fixed-sample or sequential 
data. The signal may cause a men shift and/or a change in the covariance structure. The score statistic can capture both spatial and temporal structures of the change and hence is particularly powerful in detecting weak signals. The score statistic is computationally efficient and statistically powerful. Our main theoretical contribution are  accurate analytical approximations on the false alarm rate of the detection procedures, which can be used to calibrate the threshold analytically. Numerical experiments on simulated and real data demonstrate the good performance of our procedure for solar flame detection and water quality monitoring.
\end{abstract}

\section{Introduction}

Detection the emergence of a signal in noisy background arises in many multi-sensor spatio-temporal surveillance applications. When the monitored process is in-control, sensors observe noise. When the monitored process is out of control, a signal is added to the noise, which typically possesses particular spatial and temporal correlation structure. One application is the environmental sensor network, which is used to monitor of river systems to detect a potential contaminant hazard \citep{kim2017impact}. When the signal emerges, observations from sensors may have a time-varying mean and a complicated spatio-temporal correlation structure due to water flow.

Exploiting spatio-temporal structures of the change may  enable us to detect weak signals. However, most existing methods only capture either spatial correlation \citep{healy1987note,crosier1988multivariate,jiang2011spatio-temporal,lee2014spatio-temporal,lee2015robust} or temporally correlation \citep{xie2012spectrum}. It is still not clear how to jointly capturing the spatial and temporal information in detection statistics. Moreover, computational complexity is often a concern for sensor network applications since there can be a large number of sensors. 
One issue with the classic likelihood ratio statistic is that in forming the statistics, one has to invert its sample covariance matrix, which causes both computational instability and complexity.
An alternative to the likelihood ratio statistic is the score statistic, which has also been used for developing detection procedures. When the hypothesis test is for a univariate parameter, the score test is the locally most powerful test \citep{rao1946locally}. 

We propose a new efficient score statistic for spatial-temporal surveillance, which we call the $\textsf{S}^3\textsf{T}$ statistic.  The  $\textsf{S}^3\textsf{T}$ statistic can capture both spatial and temporal correlation of a possible change signal. Hence, it can react quickly to a change in the mean and/or in the spatio-temporal covariance. An appealing feature of the score statistic is that it avoids computing the inversion of a sample covariance matrix, which leads to high computational efficiency for high-dimensional problems.  Our main theoretical contributions are accurate analytic approximations for the false detection rate in the offline case or the average run length for the online case, so calibrating thresholds to control the false alarm rate of our procedure can be done efficiently without resorting to onerous numerical simulation. This is useful in practice, as the usual trial-and-error approach to calibrate thresholds by simulation can be quite time-consuming, especially in the high-dimensional setting. When we have scalar observations (the dimension of the observation is one), our statistic $\textsf{S}^3\textsf{T}$ reduces to the score detector considered in \citep{xie2012spectrum}. In this sense, our work provides a novel and highly nontrivial extension of \citep{xie2012spectrum} when there are both spatial and temporal correlations.

The rest of the paper is organized as follows. Section 2 formulates the problem. Section 3 presents our $\textsf{S}^3\textsf{T}$ statistic for offline detection and contains theoretical approximation for the significance level and verifies its accuracy by simulations. Section 4 extends our detection procedure to the online setting and presents accurate approximation to the average-run-length. Section 5 contains numerical results that demonstrate the good performance of our procedure for solar flare detection and water quality monitoring using sensor networks. Finally, Section 6 concludes the paper.

\section{Problem Formulation} \label{formulation}

Consider a sequence of samples $\bm{y}_{\ell} \in \mathbb{R}^{p}$, $\ell=1,2,\cdots,N$, where $p$ is the dimension, $N$ is the sample size, which is fixed in the offline setting, and is growing in the online setting. We assume that under the null hypothesis, $\{\bm{y}_{\ell}\}$ forms a temporally i.i.d.\@ random noise process with spatial correlation that is caused by, for instance, either sensor measurement errors or background noises from the environment. At some time $k$, which corresponds to the {\it unknown} change-point, a signal emerges over the observation noise. The change may alter not only the mean of $\{\bm{y}_{\ell}\}$ but also the spatio-temporal correlation structure.

We first consider an offline detection setting, with the goal to detect the emergence of a change in retrospect using offline samples.  Formally, this  can be formulated as the following hypothesis test:
\begin{equation*}
\begin{array}{ll}
\textsf{H}_{0}: & \quad\  \bm{y}_{\ell} = \bm{w}_{\ell}, \qquad\quad \ell=1,2,\cdots,N, \\
\textsf{H}_{1}: & \left\{ \begin{array}{ll}
\bm{y}_{\ell} = \bm{w}_{\ell}, & \ell=1,2,\cdots,k, \\
\bm{y}_{\ell} = \bm{x}_{\ell} + \bm{w}_{\ell}, & \ell=k + 1,\cdots,N\, \end{array} \right.
\end{array}
\end{equation*}
where  $\bm{w}_{\ell} \stackrel{i.i.d.}{\sim} \mathcal{N}(0, \bm{\Sigma})$ and $\bm{\Sigma}$ is the spatial covariance matrix of the noise. Before the change, there is no temporal correlation among the samples. This is a reasonable, because usually we have plenty of data before change to estimate the temporal correlation and perform ``whitening'' to remove the temporal correlation when there is no change. 

Below we describe models for the underlying signal $\{\bm{x}_\ell\}$ when the change occurs. The signal can be spatially and temporally correlated. For temporal correlation, we use various multivariate time-series models. For instance, we  consider  the first-order vector autoregressive VAR(1) model \citep{BrockwellTimeseries}, 
\[\bm{x}_{\ell} = \bm{\mu}_{x} + \theta \bm{x}_{\ell-1} + \bm{\epsilon}_{\ell},, \ell = 1, 2, \ldots \] where $\theta \in \mathbb{R}$ and $\epsilon_{\ell} \in \mathbb{R}^{p}$ is the process noise which drives the randomness of the signal. Another example is the VARMA(1, 1) model, which is given by
 \[\bm{x}_{\ell+1} + \phi \bm{x}_{\ell}= \bm{\mu}_{x} + \eta \bm{\epsilon}_{\ell} + \bm{\epsilon}_{\ell+1},\] where $\eta \in \mathbb{R}$ and $\phi \in \mathbb{R}. $
For spatial correlation, we adopt standard spatial statistics models \citep{GaetanSpatialModel}. Denote $\E[\bm{x}_{\ell}] = \bm{\mu}_{\ell} \in \mathbb{R}^{p}$ and $\Var(\bm{x}_{\ell})=\gamma \bm{\Lambda} \in \mathbb{R}^{p \times p}$, where $\bm{\Lambda}$ is the spatial correlation matrix of the signal $\bm{x}_{\ell}$, and $\gamma \in \mathbb{R} \geq 0$ is the magnitude of the covariance of the signal. Here we assume  stationarity of the spatial covariance and  that the {\it structure} of $\bm{\Lambda}$ is known but the parameter value is unknown. This is a common practice in spatial statistic, because once a spatial correlation model is assumed,  $\bm{\Lambda}$ is specified by the location of the samples and the unknown value of parameters in the spatial model. In particular, each entry of the spatial covariance $\bm{\Lambda}$ is determined by a correlation function, 
$C(d|\rho)$, 
which is a function of the distance $d$ between two samples (in our case, sensors) and is parameterized by (unknown valued) $\rho$. Moreover, we assume the signal magnitude $\gamma$ is unknown. 
Some commonly used spatial models are as follows. Let $\mathbbm{1}\{A\}$ denotes the indicator function which takes value 1 when the indicated event $A$ is true and 0 otherwise. 
\begin{enumerate}
	
	\item Spherical model \citep{lee2014spatio-temporal}:
	\begin{equation} \label{spherical_model}
	C(d|\rho) = 1 \mathbbm{1}\{d = 0\} + \rho \mathbbm{1}\{d = 1\} + \frac{\rho}{2}\mathbbm{1}\{d = \sqrt{2}\}, \quad \rho \in [0, 1].
	\end{equation}
	
	\item Exponential model \citep{GaetanSpatialModel}:
	$$
	C(d|\rho) = 1 \mathbbm{1}\{d = 0\} + e^{-d/\rho} \mathbbm{1}\{d >0\}, \quad \rho > 0.
	$$
	
	\item Mat\'ern model \citep{GaetanSpatialModel}:
	$$
	C(d|\rho) = 1 \mathbbm{1}\{d = 0\} +  \frac{1}{2^{v-1}\Gamma(v)}(\sqrt{2}v^{1/2}d/\rho)^{v}K_{v}(\sqrt{2}v^{1/2}d/\rho)\mathbbm{1}\{d > 0\}, \quad \rho>0.
	$$
	where the parameters $\rho > 0$ and $\theta>0$, $\Gamma(\cdot)$ is the gamma function, $K_{v}(\cdot)$ is the modified Bessel function of the second kind \citep{ripley2005spatial}, $v$ is the order of the Mat\'ern model, which determines the degree of smoothness of the correlation function. Note that when $v = p + 0.5$, $p \in \mathbb{R}^+$, the Mat\'ern model can be written as a product of an exponential and a polynomial of order $p$. When $v = 0.5$,  the Mat\'ern model is equivalent to the exponential model, and when $v \rightarrow \infty$, it converges to the squared exponential covariance function.
\end{enumerate}

Now we derive our detection statistic. 
%
%
For an assumed change location $k$, let $\tau = N - k$ denote the number of post-change samples. Define a vector by stacking all post-change samples:\begin{equation} \label{y_vector}
\bm{y}_{(k+1:N)} = [\bm{y}_{k+1}^\intercal,\cdots,\bm{y}_{N}^\intercal]^\intercal \in \mathbb{R}^{p\tau},
\end{equation}
where $a^\intercal$ denotes the transpose of a vector $a$. Define $\bm{x}_{(k+1:N)}$ and $\bm{w}_{(k+1:N)}$ in a similar fashion. Then we have \[\bm{y}_{(k+1:N)} =  \bm{x}_{(k+1:N)} + \bm{w}_{(k+1:N)}.\] Now the covariance matrix of the stacked observation vector can be shown to consist of two terms that are due to the signal and the noise, respectively: \[\Var[\bm{y}_{(k+1:N)}] = \gamma \bm{V}_\tau(\theta) + \bm{\Sigma}_\tau,\] where $\bm{\Sigma}_\tau = \Var[\bm{w}_{(k+1:N)}]$, $\gamma \bm{V}_\tau(\theta) = \Var[\bm{x}_{(k+1:N)}]$ and $\theta$ is the parameter related to the temporal correlation which we will specify in the next paragraph. The second term in the covariance matrix  is given by
\begin{equation} \label{Sigma_w}
\bm{\Sigma}_\tau =
\bm{I}_\tau \otimes \bm{\Sigma} \in \mathbb{R}^{p\tau \times p\tau},
\end{equation}
where $\bm{I}_\tau$ is a $\tau$-by-$\tau$ identity matrix and $\otimes$ denotes the Kronecker product.

Through the definition in \eqref{y_vector}, the spatial and temporal correlation of the signal is jointly captured by the matrix $\bm{V}_\tau(\theta)$ and its form can be specified explicitly. For instance, for VAR(1), we have
\begin{equation} \label{V_AR1}
\bm{V}_\tau(\theta) =
\bm{R}_{\tau}(\theta) \otimes \bm{\Lambda},
\end{equation}
where $\bm{R}_{\tau}(\theta) \in \mathbb{R}^{\tau \times \tau}$ and $[\bm{R}_{\tau}(\theta)]_{i, j} = \theta^{|i-j|}, \forall i, j \in \{1,\cdots,\tau\}$.
%
Similarly, if the signal follows the VARMA(1,1) model, the matrix $\bm{V}$ can be parameterized by $\theta \triangleq (\phi, \eta)$ with the following form:
\begin{equation} \label{V_ARMA11}
\begin{aligned}
\bm{V}_\tau(\theta) &=
\bm{R}_{\tau}(\phi, \eta) \otimes \bm{\Lambda},
\end{aligned}
\end{equation}
where $\bm{R}_{\tau}(\phi, \eta) \in \mathbb{R}^{\tau \times \tau}$; $[\bm{R}_{\tau}(\phi, \eta)]_{i, j} = 1 + \eta^2 - 2\phi\eta$, if $i=j$ and $[\bm{R}_{\tau}(\phi, \eta)]_{i, j} = \phi^{|i-j|-1}(\phi-\eta)(1-\phi\eta)$, otherwise.
%
%
For the more general models, similar forms of $\bm V_\tau$ will hold, i.e., the temporal dependence of the signal is captured by some matrix $\bm{R}_{\tau}$, while the spatial dependence by $\bm{\Lambda}$, and the spatial-temporal covariance is a Kronecker product of corresponding spatial covariance and temporal covariance matrices \citep{ENV:ENV854}. 

Using the representation above, the detection problem can be reformulated as the following hypothesis test:
\begin{equation} \label{test1}
\begin{array}{lll}
\textsf{H}_{0}: &
\bm{y}_{(1:k)} \sim \mathcal{N}\big(0, \bm{\Sigma}_{k}\big), & \bm{y}_{(k+1:N)} \sim \mathcal{N}\big(0, \bm{\Sigma}_\tau\big),\\
\textsf{H}_{1}: &
\bm{y}_{(1:k)} \sim \mathcal{N}\big(0, \bm{\Sigma}_{k}\big), & \bm{y}_{(k+1:N)} \sim \mathcal{N}\big(\bm{\mu}_{(k+1:N)} ,\  \gamma \bm{V}_\tau(\theta) + \bm{\Sigma}_\tau\big), \end{array} \\
\end{equation}
where $\bm{\mu}_{(k+1:N)} = [\bm{\mu}^\intercal_{k+1},\cdots,\bm{\mu}^\intercal_{T}]^\intercal \in \mathbb{R}^{p\tau}$ and $\gamma \in \mathbb{R} >0$. Note that the above hypothesis test is equivalent to the following  simpler form that will enable us to derive the score statistic:
\begin{equation*}
\begin{aligned}
\textsf{H}_{0}:\ \  &\gamma=0,\ \bm{\mu}_{(k+1:N)} =0,\\
\textsf{H}_{1}:\ \  &\gamma >0,\ \bm{\mu}_{(k+1:N)} \neq 0,
\end{aligned}
\end{equation*}
where $\bm{a} \neq 0 $ denotes an element-wise inequality.

\section{$\textsf{S}^3\textsf{T}$ Statistic for Offline Detection} \label{procedure}

In this section we derive the $\textsf{S}^3\textsf{T}$ statistic for detection for the offline setting, i.e., all samples are collected and we aim to distinguish two hypothese.
The log-likelihood function of the hypothesis test in \eqref{test1} is given by
\begin{equation} \label{lr_fun}
\begin{aligned}
\ell(\gamma, \bm{\mu}, \tau, \theta) =& -\frac{1}{2} \log(2\pi)-\frac{1}{2} \log\big|\gamma \bm{V}_\tau(\theta) + \bm{\Sigma}_\tau\big|  \\ &- \frac{1}{2} (\bm{y}_{(k+1:N)} - \bm{\mu}_{(k+1:N)})^\intercal(\gamma \bm{V}_\tau(\theta) + \bm{\Sigma}_\tau)^{-1} (\bm{y}_{(k+1:N)} - \bm{\mu}_{(k+1:N)}).
\end{aligned}
\end{equation}
To coping with unknown parameters, one may construct a detection procedure using the generalized likelihood ratio (GLR) statistic based on (\ref{lr_fun}). However, the GLR statistic involves the calculation of the inverse of a $p\tau$-by-$p\tau$ dimensional matrix $\gamma \bm{V}_{\tau}(\theta) + \bm{\Sigma}_{\tau}$. Moreover, since the change-point location unknown, when forming the generalized likelihood ratio statistic, we have to  search over all possible  change-point locations, for $k = 1, \ldots, N$. However, for each $k$ value, calculating $(\gamma \bm{V}_{\tau}(\theta) + \bm{\Sigma}_{\tau})^{-1} $ is expensive when the dimensionality of samples $p$ or the sample size $N$ is large. Below, we show an alternative approach based on the score statistic can avoid this issue.

\subsection{Quadratic score statistic}

We now derive the score-statistic for detection. The efficient score of the model is calculated by taking the derivative of $\ell(\gamma, \bm{\mu}, \tau, \theta)$ with respect to $\gamma$ and $\bm{\mu}$ and evaluated at $\gamma=0$ and $\bm{\mu}=\bm{0}$:
\begin{equation} \label{score}
\begin{aligned}
&\varsigma(\tau, \theta) = \left[\begin{array}{c}
\frac{\partial \ell}{\partial \gamma} \big|_{\bm{\mu}=\bm{0}, \gamma=0} \\
\frac{\partial \ell}{\partial \bm{\mu}} \big|_{\bm{\mu}=\bm{0}, \gamma=0} \\
\end{array}\right] = \left[\begin{array}{c}
-\frac{1}{2}\tr\big(\bm{\Sigma}_\tau^{-1} \bm{V}_\tau(\theta)\big) + \frac{1}{2} \bm{y}_{(k+1:N)}^\intercal \bm{\Sigma}_\tau^{-1} \bm{V}_\tau(\theta) \bm{\Sigma}_\tau^{-1} \bm{y}_{(k+1:N)} \\
\bm{\Sigma}_\tau^{-1} \bm{y}_{(k+1:N)}
\end{array}\right],
\end{aligned}
\end{equation}
where $\tr\big(\cdot\big)$ denotes the trace of a matrix. The derivation of (\ref{score}) is given in the Appendix.  It can be verified that the mean of the efficient score vector $\E[\varsigma(k,\theta)]$ is $\bm{0}$ under the null hypothesis, where $\bm{0}$ represents a vector of zeros. 
It can be shown that the covariance of the score vector $\varsigma(\tau, \theta)$ is given by
$$
\mbox{Cov}[\varsigma(\tau, \theta)] =
 \left[\begin{array}{cc}
\frac{1}{2} \tr\Big(\bm{\Sigma}_\tau^{-1} \bm{V}_\tau(\theta) \bm{\Sigma}_\tau^{-1} \bm{V}_\tau(\theta)\Big) & 0
\\
\bm{0} & \bm{\Sigma}_\tau^{-1}
\end{array}\right].
$$
As suggested by the seminal work of \cite{radhakrishna_rao_1948}, 
when the likelihood function depends on multiple parameters, the score statistic corresponds to  a  quadratic function of the efficient score vector. In our case, this corresponds to 
\begin{equation} \label{W1}
\begin{aligned}
S(\tau, \theta) &= \varsigma(\tau, \theta)' \mbox{Cov}[\varsigma(\tau, \theta)]^{-1} \varsigma(\tau, \theta)\\
&= \frac{\Big[\bm{y}_{(k+1:N)}^\intercal \bm{\Sigma}_\tau^{-1} \bm{V}_\tau(\theta) \bm{\Sigma}_\tau^{-1} \bm{y}_{(k+1:N)} - c(\tau, \theta) \Big]^{2}}{d(\tau, \theta)} + \bm{y}_{(k+1:N)}^\intercal \bm{\Sigma}_\tau^{-1} \bm{y}_{(k+1:N)},
\end{aligned}
\end{equation}
where
$
c(\tau, \theta) = \tr\big(  \bm{\Sigma}_\tau^{-1} \bm{V}_\tau(\theta)\big),
$
and
$
d(\tau, \theta) = 2 \tr\big[\bm{\Sigma}_\tau^{-1} \bm{V}_\tau(\theta) \bm{\Sigma}_\tau^{-1} \bm{V}_\tau(\theta)\big].
$
Note that the computation of $S(\tau, \theta)$ is relatively easy and much less expensive than the GLR statistic. The only place requires matrix inversion is the inversion of $\bm{\Sigma}_\tau$. The matrix  $\bm{\Sigma}_\tau$ defined in (\ref{Sigma_w}) has a simple block diagonal structure $\bm{\Sigma}_\tau = \bm{I}_\tau \otimes \bm{\Sigma} $. Hence, its inversion only requires to compute $\bm{\Sigma} ^{-1}$, with a computational complexity of $\mathcal O (p^3)$ (which is much smaller than $\mathcal O (p\tau)^3$, if we have to directly invert $\bm{\Sigma}_\tau$). Moreover, since $\Sigma$ is assumed known and fixed, its inversion can be pre-computed and does not cause an issue for the online computation of the detection statistic. 

Since $S(\tau, \theta)$ has an increasing mean as the change location $\tau$ decreases, it needs to be normalized to have mean 0 and variance 1 under the null hypothesis. We refer to the resulting detection statistic  as the {\it quadratic score statistic}, 
\begin{equation} \label{W1_tilde}
\begin{aligned}
\widetilde{S}(\tau, \theta)  = \frac{S(\tau, \theta) - \E\big[S(\tau, \theta)\big]}{\sqrt{\Var \big[S(\tau, \theta) \big]}}.
\end{aligned}
\end{equation}
It can be shown that the mean is given by
$
\E\big[S(\tau, \theta)\big] = p\tau + 1,
$
and the variance is given by
\begin{equation*}
\begin{aligned}
\Var \big[S(\tau, \theta) \big] =& 2p\tau + 10 - 24 \frac{c(\tau, \theta)}{d(\tau, \theta)^2} \tr \Big( \bm{\Sigma}_\tau^{-1} \bm{V}_\tau(\theta)\bm{\Sigma}_\tau^{-1} \bm{V}_\tau(\theta)\bm{\Sigma}_\tau^{-1} \bm{V}_\tau(\theta) \Big) \\
&+ \frac{48}{d(\tau, \theta)^2} \tr \Big( \bm{\Sigma}_\tau^{-1} \bm{V}_\tau(\theta)\bm{\Sigma}_\tau^{-1} \bm{V}_\tau(\theta) \bm{\Sigma}_\tau^{-1} \bm{V}_\tau(\theta)\bm{\Sigma}_\tau^{-1} \bm{V}_\tau(\theta) \Big).
\end{aligned}
\end{equation*}
Then we may construct an offline detection procedure using $\widetilde{S}(\tau, \theta)$, which detects a signal when the maximum standardized score statistic over all possible $\theta$ and $\tau$ exceeds a pre-specified threshold $b$:
\[ \max_{\theta \in \Theta,\ 1 \leq \tau \leq N} \widetilde{S}(\tau, \theta) \geq b,\] where $\Theta$ is the set of possible values of the parameter $\theta$.

\subsection{$\textsf{S}^3\textsf{T}$ statistic for offline change-point detection}

Although it is claimed by \citep{radhakrishna_rao_1948} that the quadratic score statistic achieves the maximum discrimination between the null and the alternative, the statistic is too complicated to 
perform theoretical analysis and difficult to calibrate the threshold $b$. In this section, we propose a simpler statistic, namely $\textsf{S}^3\textsf{T}$ statistic, which is the score with respect to $\gamma$ {\it only}, which is given by
\begin{equation} \label{W2}
\begin{aligned}
W&(\tau, \theta) =\frac{\frac{\partial \ell}{\partial \gamma} \big|_{\bm{\mu}=\bm{0}, \gamma=0} }{\sqrt{ \Var\big[\frac{\partial \ell}{\partial \gamma} \big|_{\bm{\mu}=\bm{0}, \gamma=0}\big] }}  = \frac{\bm{y}_{(N - \tau +1:N)}^\intercal \bm{\Sigma}_\tau^{-1} \bm{V}_\tau(\theta) \bm{\Sigma}_\tau^{-1} \bm{y}_{(N-\tau+1:N)} - c(\tau, \theta)}{\sqrt{d(\tau, \theta)}}.
\end{aligned}
\end{equation}
Under the null hypothesis, the detection statistic $W(\tau, \theta)$ has mean 0, and variance 1.
Similarly, the procedure claims to detect a signal if the maximum score statistic exceed a pre-specified threshold $b$,
\begin{equation} \label{maxW2}
\begin{aligned}
\max_{\theta \in \Theta,\ 1 \leq \tau \leq N} W(\tau, \theta) \geq b.
\end{aligned}
\end{equation}

\subsection{Control false alarms of offline statistic}

In this section, we present a theoretical approximation for the significance level of the detection procedure defined in \eqref{maxW2}, which avoid the time-consuming simulation when deciding the appropriate $b$.

Below, let $\bm{A}_{\tau}(\theta) = \bm{\Sigma}_{\tau}^{-1} \bm{V}_{\tau}(\theta)$, and $\bm{B}_{\tau}(\theta) = \bm{\Sigma}_{\tau}^{-1/2} \bm{V}_{\tau}(\theta)\bm{\Sigma}_{\tau}^{-1/2}$. Let $\bm{I}_p$ denote a $p$ by $p$ identity matrix. Denote the standard normal density function by $\phi(x)$ and its distribution function by $\Phi(x)$, and define a special function \citep{Siegmund2007}: 
\begin{equation} \label{nu_fun}
\nu(x) = \frac{\frac{2}{x} \Big[\Phi\big(\frac{x}{2}\big) - \frac{1}{2}\Big]}{ \frac{x}{2} \Phi \big(\frac{x}{2}\big) + \phi\big(\frac{x}{2}\big) }.
\end{equation}
Define the following quantities, which are useful for the statement of the theorem
	\begin{align}
	\mu(\tau, \theta) &= \tau\Bigg[ \frac{\tr \big(\bm{A}_{\tau+1}(\theta) \bm{A}_{\tau+1}(\theta)\big)}{\tr \big(\bm{A}_{\tau}(\theta) \bm{A}_{\tau}(\theta)\big)} - 1 \Bigg], 	 \label{mu_fun} \\
	H(\tau, \theta)  &= - \frac{\partial^2 E[W(\tau, \theta) W(\tau, s)]}{\partial^2 s} \bigg|_{s = \theta}, 	\label{H} \\
	g(\tau, \theta) &= \frac{\exp \big(-\xi_{0}(\tau, \theta)b + \psi(\xi_{0}(\tau, \theta)  \big)}{  \sigma_{\xi_{0}} \sqrt{2\pi} },  \label{g_fun} \\
	\psi(\xi) & = -\xi \frac{c(\tau, \theta)}{ \sqrt{d(\tau, \theta)}} - \frac{1}{2} \log \bigg| \bm{I}_{p\tau} - \frac{2\xi \bm{B}_{\tau}(\theta) }{ \sqrt{d(\tau, \theta)} } \bigg|. 	 \label{psi_fun}
	\end{align}
	Note that $\psi(\xi)$ is the cumulant generating function of the detection statistic $W(\tau,\theta)$. The following theorem is one of the main theoretical contribution, which provides an analytical approximation for significance level of the detection procedure defined in \eqref{maxW2}.
\begin{theorem}[Approximation for significance level] \label{theorem1}
	When the threshold $b \rightarrow \infty$ and $\theta \in \Theta \subset \mathbb{R}^{d}$, under the null hypothesis, the probability of false detection  for the procedure defined in \eqref{maxW2} is given by
	\begin{equation} \label{SL}
	\begin{aligned}
	& \mathbb{P}_{\textsf{H}_{0}}\Big( \max_{\substack{\theta \in \Theta \\ 1 \leq \tau \leq N}} W(\tau, \theta) \geq b \Big) \\&= \frac{1}{(2\pi)^{\frac{d}{2}}} \sum_{\tau =1}^{N} \int_{\theta \in \Theta} \frac{[b \xi_{0}(\tau, \theta)]^{\frac{d}{2}}}{\xi_{0}(\tau, \theta)}g(\tau, \theta) |H(\tau, \theta)|^{\frac{1}{2}} \frac{b^{2}\mu(\tau, \theta)}{2\tau} \nu \Big(\sqrt{\frac{b^{2} \mu(\tau, \theta)}{\tau}}\Big) d\theta + o(1),
	\end{aligned}
	\end{equation}
		where
	\begin{align*}
	\sigma^2_{\xi_{0}}
	&=  d(\tau, \theta)^{-1} \tr\Big( \bigg[\bm{I}_{p\tau} - \frac{2\xi_{0} \bm{B}_{\tau}(\theta) }{\sqrt{d(\tau, \theta)} }\bigg]^{-1} \bm{B}_{\tau}(\theta) \bigg[\bm{I}_{p\tau} - \frac{2\xi_{0} \bm{B}_{\tau}(\theta) }{\sqrt{d(\tau, \theta)}}\bigg]^{-1} \bm{B}_{\tau}(\theta)\Big),
	\end{align*}
	and $\xi_{0} (\tau, \theta)$ is the solution to
	\begin{equation}	\label{equation}
	\frac{1}{\sqrt{d(\tau, \theta)}}\tr\Big( \Big[\bm{I}_{p\tau} - \frac{2\xi_{0} \bm{B}_{\tau}(\theta) }{ \sqrt{d(\tau, \theta)} }\Big]^{-1} \bm{B}_{\tau}(\theta) - \bm{A}_{\tau}(\theta) \Big) = b.
	\end{equation}
	
\end{theorem}

 The main proof technique for Theorem \ref{theorem1} involves the change-of-measure to evaluating the boundary hitting probability of random process \citep{siegmund2013sequential,yakir2013extremes}. See the Appendix for the derivation of (\ref{psi_fun}) and the proof of Theorem \ref{theorem1}, when the dimension of parameter $\theta$ is 1 (i.e.\@, $d=1$), and the proof for one-dimensional parameter space can be generalized to the multi-dimensional case using similar proof techniques.



We verify the accuracy of the approximation in Theorem \ref{theorem1} by comparing the approximated significance levels with simulated ones. In the experiment, we assume that the temporal correlation structure of the signal $\{ \bm{x}_{\ell} \}$ follows a VAR(1) model, $\bm{x}_{\ell} = \bm{\mu}_{x} + \theta \bm{x}_{\ell-1} + \bm{\epsilon}_{\ell}$,  where $\theta \in \mathbb{R}$, which means $\bm{V}_\tau(\theta)$ has the form in \eqref{V_AR1}. We further assume the spatial correlation of the signal follows a spherical model, as defined in \eqref{spherical_model}, with parameter  $\rho=0.3$.
We set $N = 50$. The search space of $\theta$ is set as a uniform grid from 0.1 to 0.9 with a step size 0.1. We vary the dimension of the signal $p$. In addition, the covariance matrix of the noise process $\bm{\Sigma}$ is assumed to be a $p$-by-$p$ identity matrix. Simulation results are based on 5000 independent replications. Both simulated and approximated false alarm rates are reported in Table \ref{tab:SL_AR1}. As one can observe, the approximation is quite accurate.

\begin{table*}[h!]
	\centering
	\caption{Simulated and approximated significance level when the signal $\{\bm{x}_{\ell}\}$ follows a VAR(1) model ($\theta \in [0.1, 0.9]$, $N = 50$ and $\rho=0.3$).}
	\begin{tabular}{r|rr|rr|rr}
		\toprule
		& \multicolumn{2}{c}{$p =2$} & \multicolumn{2}{c}{$p =9$} & \multicolumn{2}{c}{$p =36$} \\
		\midrule
		$b$     & Simulated & Approx. & Simulated & Approx. & Simulated & Approx. \\
		\midrule
		3.5   & 0.097 & 0.097 & 0.065 & 0.057 & 0.036 & 0.042 \\
		4     & 0.063 & 0.068 & 0.036 & 0.030 & 0.013 & 0.019 \\
		4.5   & 0.038 & 0.047 & 0.018 & 0.019 & 0.006 &0.008  \\
		5     & 0.033 & 0.032 & 0.011 & 0.012 & 0.003 & 0.003 \\
		5.5   & 0.022 & 0.021 & 0.005 & 0.007 & 0.002 &0.001\\
		6     & 0.015 & 0.014 & 0.003 & 0.004 & 0.0004 & 0.0005 \\
		6.5   & 0.006 & 0.009 & 0.002 & 0.002 & 0.0002 & 0.0002 \\
		\bottomrule
	\end{tabular}%
	\label{tab:SL_AR1}%
\end{table*}%


\section{$\textsf{S}^3\textsf{T}$  for online change-point detection}

In this section, we present an online change-point detection procedure based on the $\textsf{S}^3\textsf{T}$ statistic.
In the online detection setting, the total sample size $N$ is not fixed, and new observations are sequentially collected. A signal may occur at some unknown change-point $k$, and our goal is to detect the emergence of the signal as soon as it occurs. The model can be described as,

\begin{equation*}
\begin{array}{ll}
\textsf{H}_{0}: & \quad\  \bm{y}_{\ell} = \bm{w}_{\ell}, \qquad\quad \ell=1,2,\cdots, \\
\textsf{H}_{1}: & \left\{ \begin{array}{ll}
\bm{y}_{\ell} = \bm{w}_{\ell}, & \ell=1,2,\cdots,k, \\
\bm{y}_{\ell} = \bm{x}_{\ell} + \bm{w}_{\ell}, & \ell=k + 1,\cdots, \end{array} \right.
\end{array}
\end{equation*}

Here we adopt a sliding window approach for the online detection procedure. We construct detection statistics using the most recent $\omega$ samples at each time, where $\omega$ is a pre-specified window length (demonstrated in Figure \ref{sliding_window} in the appendix). We did not search for the unknown change-point location to reduce computational complexity (since this has to be done for each time $t$). This corresponds to a type of Shewhart chart \citep{Shewart1931}.  
Given a current time $t$, the detection statistic constructed using the most recent $\omega$ samples is given by
\begin{equation} \label{W_online}
\begin{aligned}
&W_{t}(\omega, \theta) = \frac{\bm{y}_{(t-\omega+1:t)}^\intercal \bm{\Sigma}_\omega^{-1} \bm{V}_\omega(\theta) \bm{\Sigma}_\omega^{-1} \bm{y}_{(t-\omega+1:t)} - c(\omega, \theta)}{\sqrt{d(\omega, \theta)}}.
\end{aligned}
\end{equation}


The detection procedure is a stopping time, which raises an alarm when the detection statistic exceeds a threshold for the first time: 
\begin{equation} \label{stopping_rule}
\mathcal{T} = \inf \bigg\{ t: \max_{\substack{\theta \in \Theta}} W_{t}(\omega, \theta) \geq b \bigg\},
\end{equation}
where $b$ is a pre-specified threshold.

\subsection{Control false alarm rate for online statistic}

In the online detection setting, the performance metric used for characterizing false alarm rate is the average-run-length (ARL), which is the expected stopping time of the procedure when there is no signal,  denoted as $E_{\textsf{H}_{0}}(\mathcal{T})$. The following theorem provides an approximation on $E_{\textsf{H}_{0}}(\mathcal{T})$ of the detection procedure defined in \eqref{stopping_rule}.

\begin{theorem}[Approximation of average-run-length] \label{theoremARL}
	Assume that $b \rightarrow \infty$. For the stopping rule defined in \eqref{stopping_rule}, the approximation on $E_{\textsf{H}_{0}}(\mathcal{T})$ is given by
	\begin{equation} \label{ARL}
	\begin{aligned}
	&E_{\textsf{H}_{0}}(\mathcal{T}) = (2\pi)^{\frac{d}{2}} \Bigg( \int_{\theta \in \Theta} \frac{[b \xi_{0}(\omega, \theta)]^{\frac{d}{2}}}{\xi_{0}(\omega, \theta)}g(\omega, \theta) |H(\omega, \theta)|^{\frac{1}{2}} \frac{b^{2}\mu(\omega, \theta)}{2\omega} \nu \Big(\sqrt{\frac{b^{2} \mu(\omega, \theta)}{\omega}}\Big) d\theta \Bigg)^{-1} \Big(1 + o(1)\Big).
	\end{aligned}
	\end{equation}
\end{theorem}

The derivation of Theorem \ref{theoremARL} uses the similar technique as in the derivation of Theorem \ref{theorem1}. By Theorem \ref{theorem1}, we can first obtain an approximation of the probability $\mathbb{P}_{\textsf{H}_{0}} ( \mathcal{T} \leq m )$, where $m$ is fixed and sufficiently large:
\begin{equation} \label{sl_online}
\begin{aligned}
& \mathbb{P}_{\textsf{H}_{0}} \Big( \mathcal{T} \leq m \Big) = \mathbb{P}_{\textsf{H}_{0}} \Big( \max_{\substack{\theta \in \Theta \\ 1 \leq t \leq m}} W_{t}(\omega, \theta) \geq b \Big) \\
&= (2\pi)^{-\frac{d}{2}} \Bigg( \sum_{t =1}^{m} \int_{\theta \in \Theta} \frac{[b \xi_{0}(\omega, \theta)]^{\frac{d}{2}}}{\xi_{0}(\omega, \theta)}g(\omega, \theta) |H(\omega, \theta)|^{\frac{1}{2}} \frac{b^{2}\mu(\omega, \theta)}{2\omega} \nu \Big(\sqrt{\frac{b^{2} \mu(\omega, \theta)}{\omega}}\Big) d\theta \Bigg) + o(1).
\end{aligned}
\end{equation}

As argued in \cite{siegmund1995using} and \cite{siegmund2008detecting}, the stopping time $\mathcal{T}$ is asymptotically exponentially distributed and is uniformly integrable. Hence, for large $m$, $\mathbb{P}_{\textsf{H}_{0}} ( \mathcal{T} \leq m ) - [1 - \exp(-\lambda m)] \rightarrow 0$, where $\lambda$ is equal to the right hand side of \eqref{sl_online} divided by $m$. Thus $E_{\textsf{H}_{0}}(\mathcal{T}) \approx \lambda^{-1}$, which is equivalent to \eqref{ARL}.

\begin{figure}[t]
	\begin{center}
		\begin{tabular}{ccc}
			\includegraphics[width=.32\textwidth]{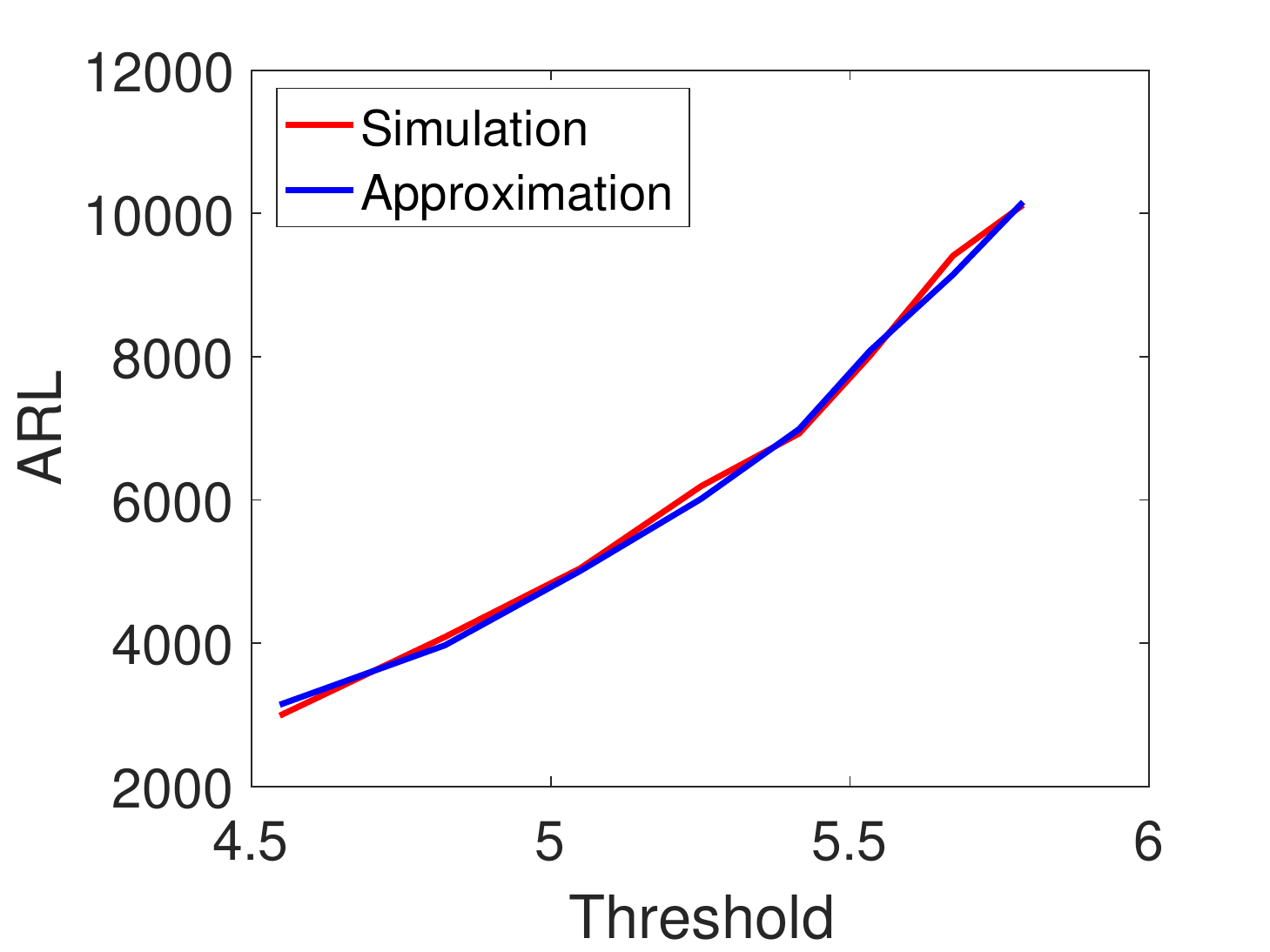} &
			\includegraphics[width=.32\textwidth]{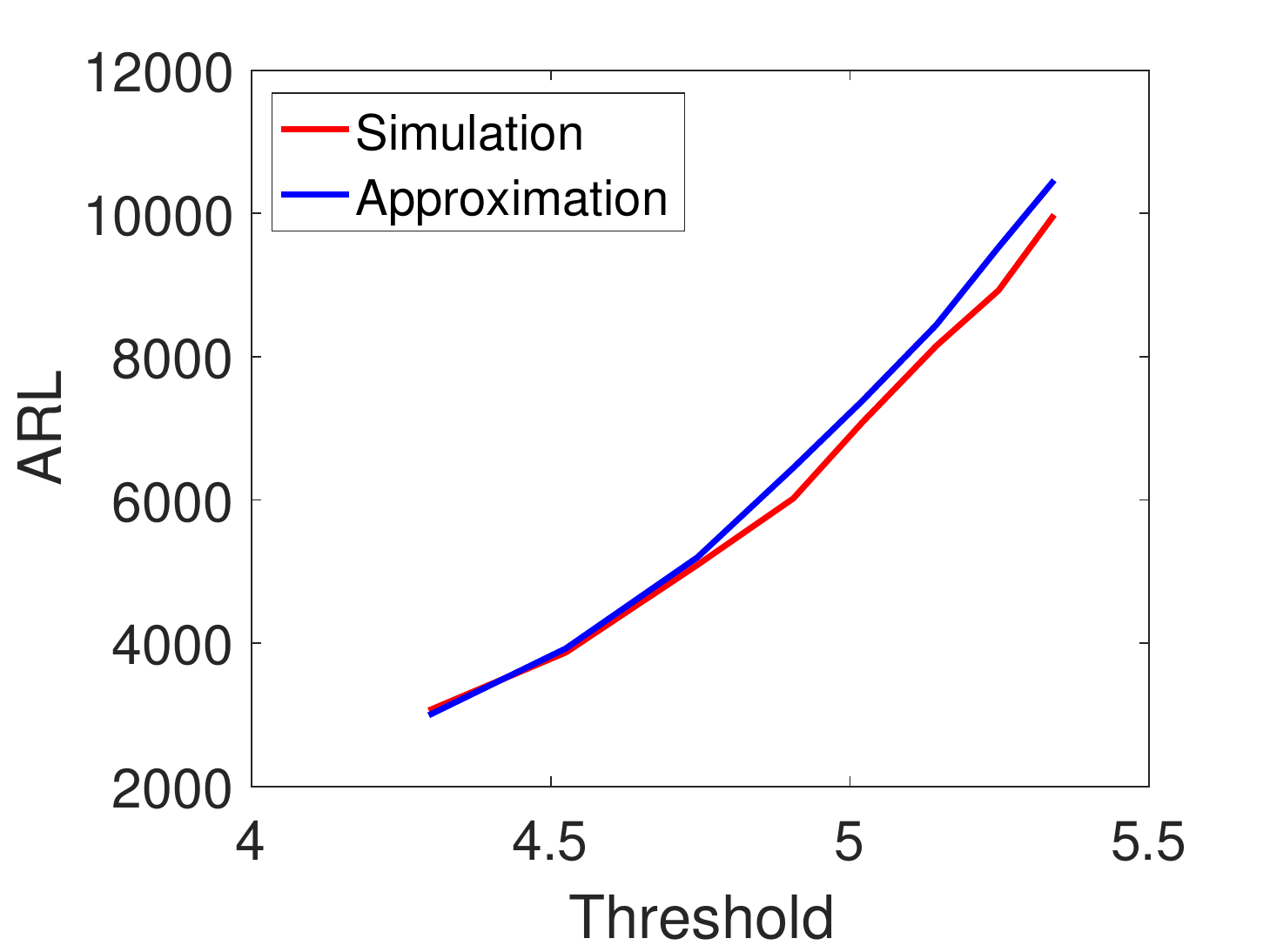} &
			\includegraphics[width=.32\textwidth]{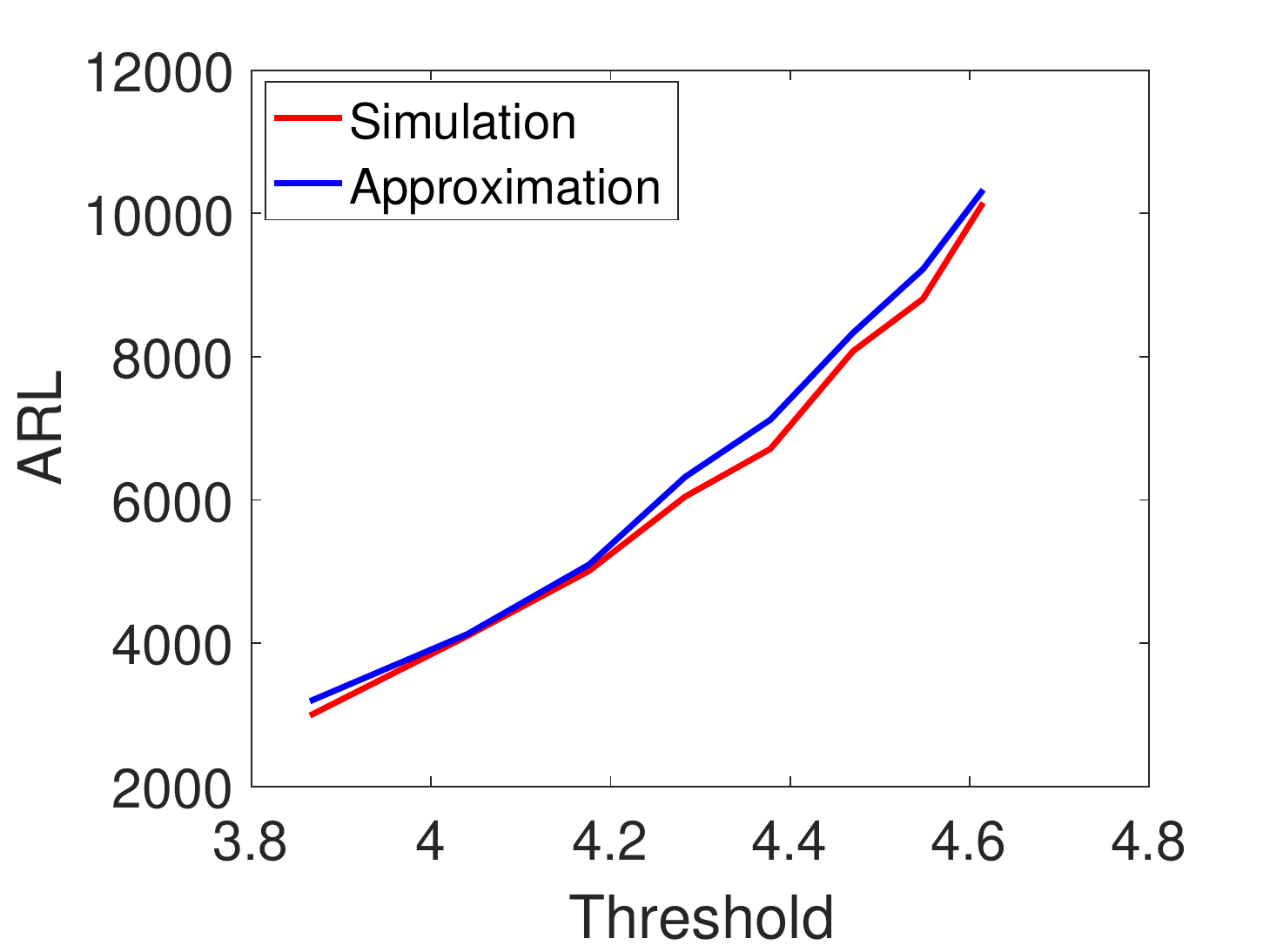} \\
			(a) & (b) & (c)
		\end{tabular}	
	\end{center}
	\caption{Comparison of approximated and simulated ARL for (a) $p=1$, (b) $p=2$, and (c) $p=9$.}
	\label{ARL_comparison}
\end{figure}

The accuracy of Theorem \ref{theoremARL} is verified by comparing simulated and approximated $E_{\textsf{H}_{0}}(\mathcal{T})$. In the experiments, the signal $\{ \bm{x}_{\ell} \}$ is generated by a VAR(1) model, $\bm{x}_{\ell} = \bm{\mu}_{x} + \theta \bm{x}_{\ell-1} + \bm{\epsilon}_{\ell}$,  where $\theta \in \mathbb{R}$. Hence, $\bm{V}_\tau(\theta)$ has the form in \eqref{V_AR1}. Meanwhile, we assume that the spatial correlation of the signal follows a spherical model, as defined in \eqref{spherical_model}, with parameter $\rho=0.3$. The search space of parameter $\theta$ is a uniform grid from 0.1 to 0.9 with interval 0.1. In addition, the covariance matrix of the noise process $\bm{\Sigma}$ is assumed to be a $p$-by-$p$ identity matrix. The results based on 5000 replications are presented in Figure \ref{ARL_comparison}.  The comparison between simulated and approximated ARLs shows that the approximation in Theory \ref{theoremARL} is quite accurate.

\section{Numerical Examples and Power Study}

In this section, we demonstrate the performance of the proposed detection statistic by comparing with other methods on simulated data, on a real-data example of solar flare detection, and on a synthetic example simulated from realistic setting for water quality monitoring.

We  focus on performance comparison for online change-point detection, since it is the most relevant setting for our targeted application. The performance comparison for offline change-point detection will be similar. We adopt the commonly used performance metric for online setting, the expected detection delay (EDD) after a change has occured. There is a tradeoff in the average run length (ARL) when there is no change and the EDD. Typically, we choose the threshold for each procedure so that its ARL meets a pre-specified large value (e.g., 5000 or 10000), so that there is rarely a false alarm. 


\subsection{Simulation}

The detection procedure defined in \eqref{stopping_rule} is compared with two other procedures: (i) an online detection procedure defined in a similar way as \eqref{stopping_rule} using the quadratic score statistic $\widetilde{S}(\tau, \theta)$, and (ii) a multivariate cumulative sum (MCUSUM) procedure \citep{healy1987note}. In the MCUSUM procedure, at each time step, a $T^2$ type of statistic \citep{hotelling1947multivariate} is calculated, and a CUSUM procedure is constructed based on the $T^2$ statistic.

In the experiment, the signal is generated from a VAR(1) model, $\bm{x}_{\ell} = \bm{\mu}_{x} + \theta \bm{x}_{\ell-1} + \bm{\epsilon}_{\ell}$, with sample dimensionality $p=2$ and parameter $\theta = 0.5$. The spatial model of the signal follows the spherical model defined in \eqref{spherical_model} with $\rho = 0.3$. For the two procedures based on $\textsf{S}^3 \textsf{T}$ and the quadratic score statistic, we use window length $\omega=50$ and search space for the parameter $\theta$, $\{0.1, 0.2, \cdots, 0.9\}$. 
Thresholds for all three procedures are calibrated so that $E_{\textsf{H}_{0}}(\mathcal{T}) = 100$. The change occurs at $t=1$.  We keep the mean of the signal $\bm{\mu}_{x} = \E[\bm{x}_{\ell}]$ as a constant (not time-varying) vector with all elements equal to $\mu$. We explore different values of $\mu$ for the mean shift and $\gamma$ for the magnitude of covariance matrix of the signal. If $\mu = 0$ and $\gamma > 0$, then the signal only causes change in covariance; if both $\mu$ and $\gamma$ are positive, then there are both mean shift and covariance change. Hence, the experiments demonstrate that the proposed detection procedure is suitable for cases where there is either mean shift or covariance change, or both.

Table \ref{tab:EDD} reports the simulated EDD of  three procedures based on 5000 replications. The smallest EDD values for each setting are marked  bold. The comparison shows that the two score type of procedures which capture both spatial and temporal correlation, i.e.\@,  $\textsf{S}^3 \textsf{T}$ and the quadratic score statistic outperform the MCUSUM procedure( which only captures the spatial correlation information). Such advantage is more significant when the signal is weak, i.e., when $\gamma$ or $\mu$ are small. This demonstrates that incorporating temporal correlation information indeed improves  detection performance. We also find that $\textsf{S}^3 \textsf{T}$  outperforms the quadratic score statistic in many cases. Given that $\textsf{S}^3 \textsf{T}$ enjoys tractable theoretical analysis and an accurate approximation for its false alarm rate, it is a good option for practitioners.

\begin{table}[h!]
	\centering
	\caption{Simulated expected detection delay.}
	\scalebox{0.8}{
		\begin{tabular}{r|rrrrr|rrrrr|rrrrr}
			\toprule
			\multicolumn{6}{c}{$\textsf{S}^3 \textsf{T}$}                         & \multicolumn{5}{c}{Quadratic score statistic}                 & \multicolumn{5}{c}{MCUSUM} \\
			\midrule
			$\gamma \backslash \mu$ & 0     & 0.1   & 0.5   & 1     & 2     & 0     & 0.1   & 0.5   & 1     & 2     & 0     & 0.1   & 0.5   & 1     & 2 \\
			\midrule
			0.01  & \textbf{97.27} & \textbf{59.08} & \textbf{6.37} & 2.80  & \textbf{1.49} & 98.05 & 65.82 & 6.45  & \textbf{2.77} & 1.51  & 98.37 & 77.67 & 9.43  & 3.56  & 1.79 \\
			0.05  & 96.28 & \textbf{57.96} & \textbf{5.95} & \textbf{2.72} & \textbf{1.49} & \textbf{95.32} & 63.19 & 6.74  & 2.81  & 1.52  & 96.79 & 71.97 & 9.28  & 3.54  & 1.79 \\
			0.1   & \textbf{72.93} & \textbf{53.16} & \textbf{6.04} & \textbf{2.78} & 1.50  & 82.49 & 56.78 & 6.74  & 2.86  & \textbf{1.49} & 80.70 & 65.16 & 9.21  & 3.54  & 1.78 \\
			0.2   & \textbf{65.32} & \textbf{46.16} & \textbf{5.96} & \textbf{2.77} & 1.50  & 74.87 & 48.83 & 6.28  & 2.78  & \textbf{1.47} & 67.33 & 55.17 & 9.02  & 3.52  & 1.79 \\
			0.5   & 39.40 & \textbf{30.32} & \textbf{5.81} & \textbf{2.78} & 1.56  & \textbf{37.07} & 33.42 & 6.07  & 2.80  & \textbf{1.50} & 41.52 & 35.87 & 8.36  & 3.47  & 1.78 \\
			1     & \textbf{20.91} & \textbf{19.42} & \textbf{5.65} & \textbf{2.75} & \textbf{1.51} & 22.75 & 20.51 & 5.64  & 2.76  & 1.55  & 23.71 & 21.31 & 7.45  & 3.45  & 1.77 \\
			\bottomrule
	\end{tabular}}%
	\label{tab:EDD}%
\end{table}%

\subsection{Solar flare detection}

\begin{figure}[h!]
	\begin{center}
		\begin{tabular}{cc}
			\includegraphics[width=.32\textwidth]{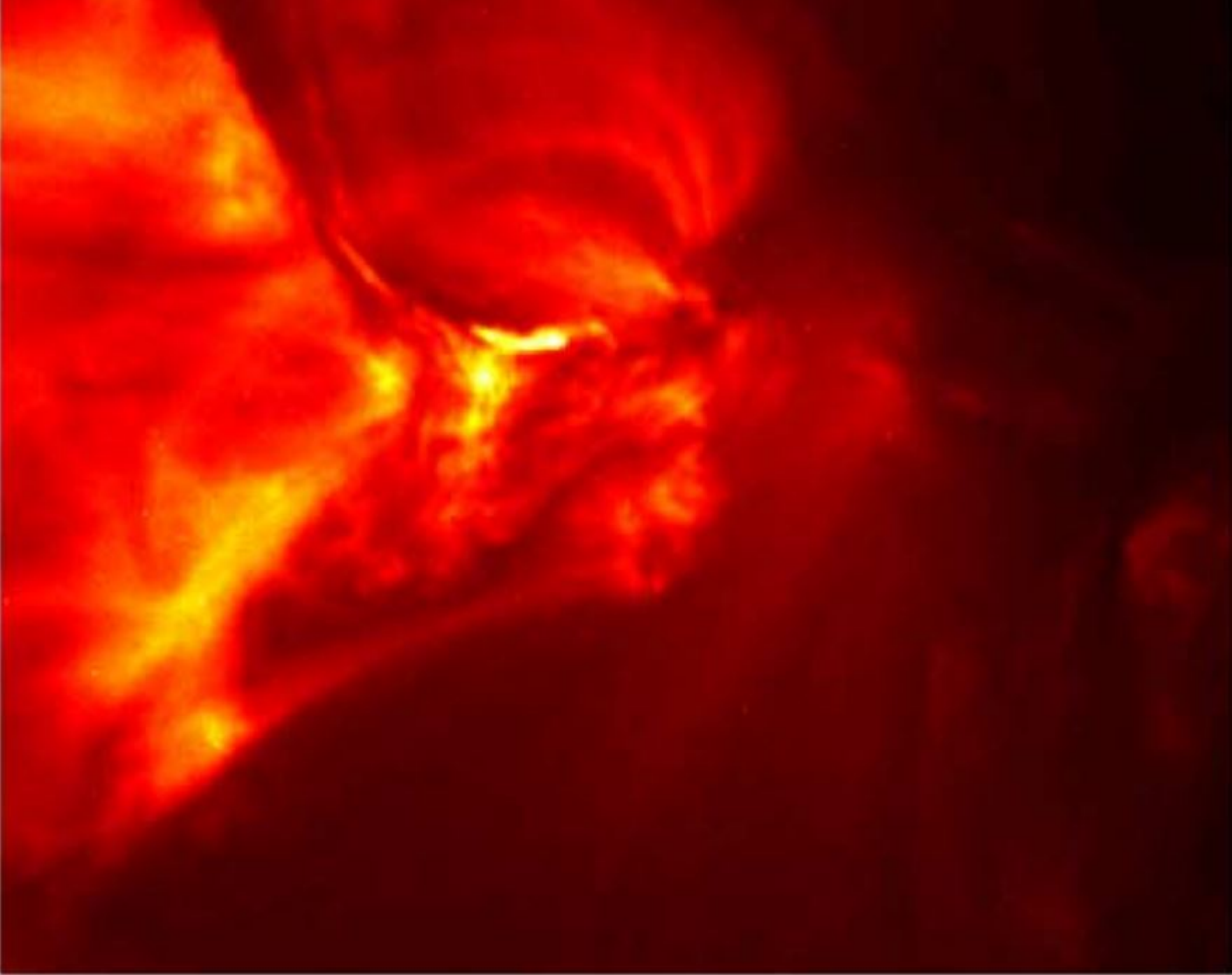} &
			\includegraphics[width=.32\textwidth]{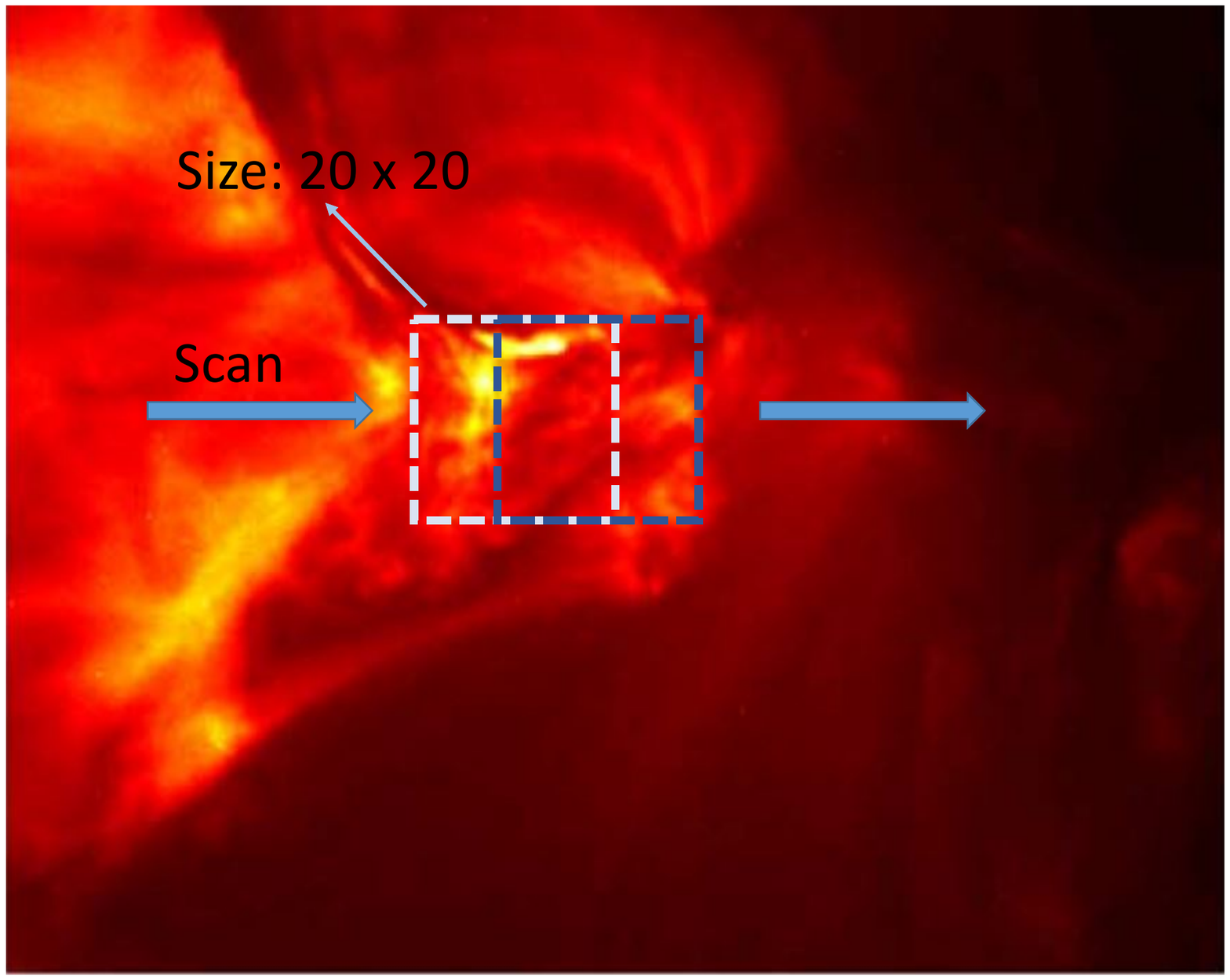}
		\end{tabular}	
	\end{center}
	\caption{Detection of solar flare at $t = 227$: (left) snapshot of the original SDO data at $t = 227$; (right) overlapping image patches for dimensionality reduction. }
	\label{solar_scan}
\end{figure}

We apply our detection procedure on a dataset from the Solar Data Observatory \citep{SDO}. The data is a video sequence that contains an abrupt emergence of a solar flare occurs around time $t = 227$. In this video, the normal states are slowly drifting image of sun surface, and the anomaly is a much brighter transient solar flare. A snapshot from this dataset during a solar flare at $t=227$ is shown in the Figure \ref{solar_scan}.

The size of the images is $\times$ 292 pixels. After vectoring the image this leads to $67744$ dimensional vectors. Due to high dimensionality, It is computationally expensive to directly apply our detection procedure on the original images. Hence, we apply spatial scanning by breaking the original image into overlapping patches of dimension $20 \times 20$, as demonstrated in the right figure of Figure \ref{solar_scan} in the appendix. The detection statistic is calculated for each image patch (of dimension $p = 400$), and we take the maximum among all patches as the detection statistic. 

\begin{figure}[t]
	\begin{center}
		\begin{tabular}{cc}
			\includegraphics[width=.32\textwidth]{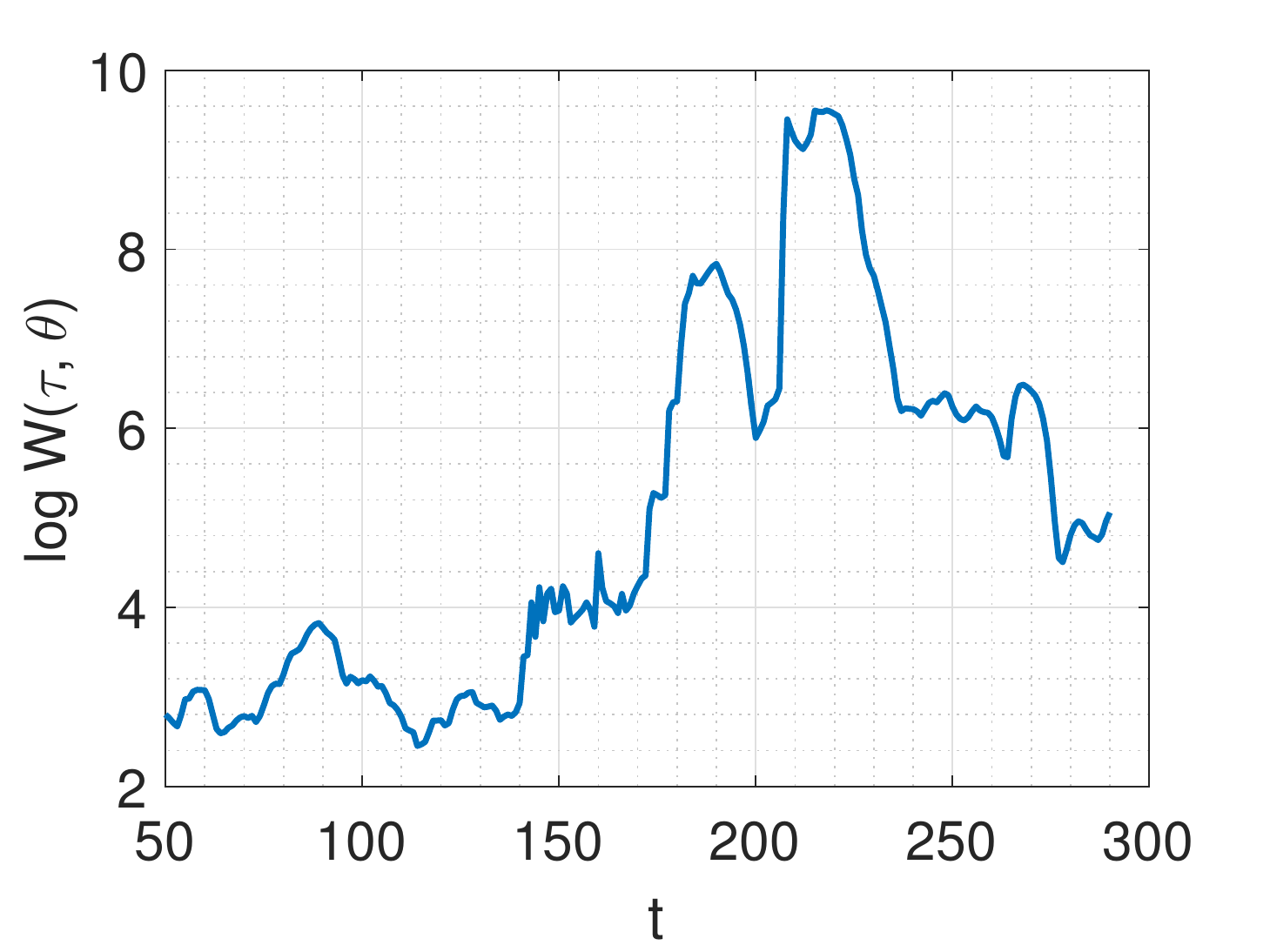} &
			\includegraphics[width=.32\textwidth]{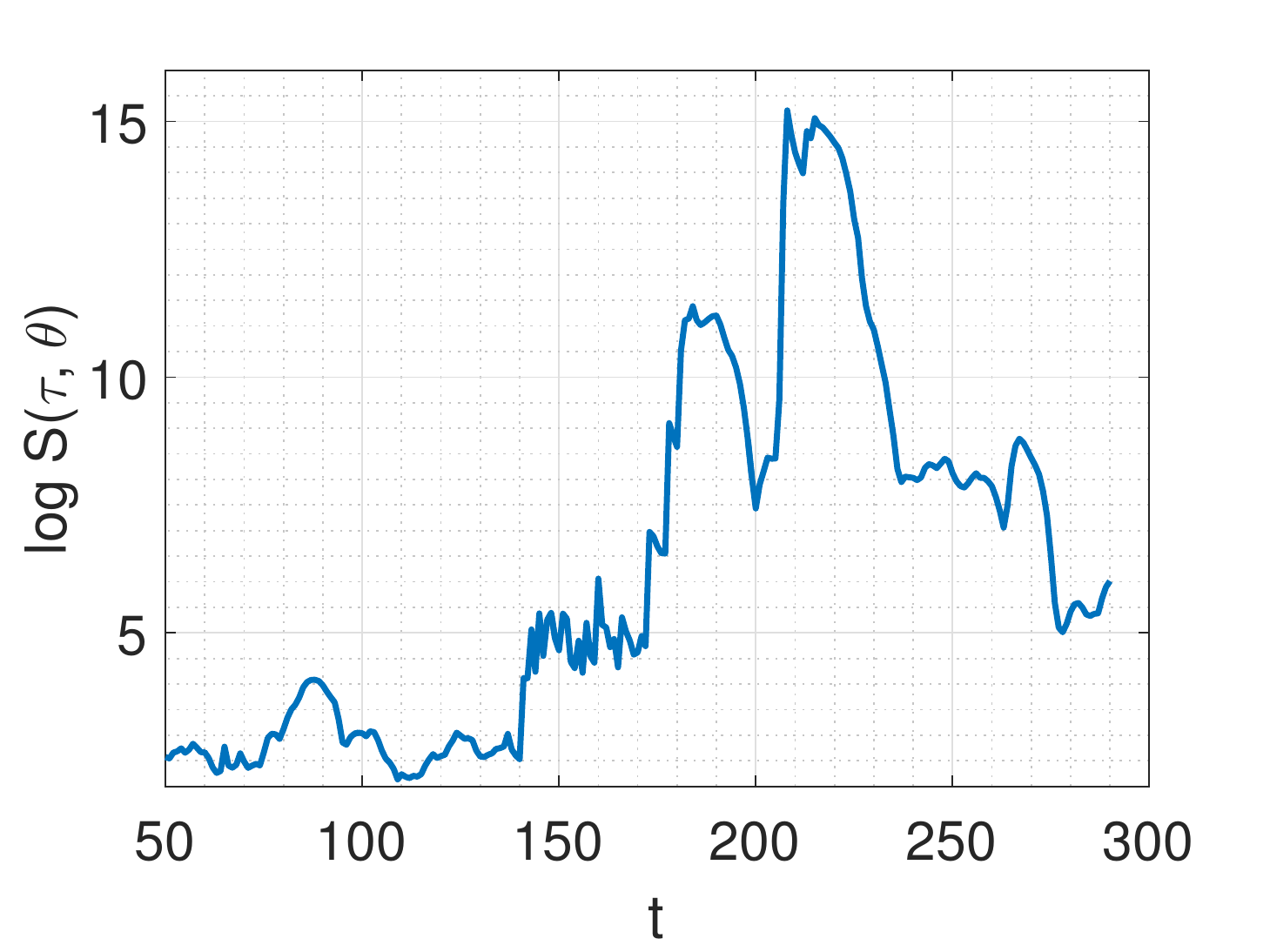}\\
			(a) $\textsf{S}^3 \textsf{T}$ & (b) Quadratic score statistic
		\end{tabular}	
	\end{center}
	\vspace{-0.2in}
	\caption{Detection statistics on logarithmic scale.}
	\label{solar_stat}
\end{figure}

We assume that before the solar flare, the data form a white noise process with no spatial and temporal correlation. The mean and variance of the noise process are estimated by the first 50 samples in the sequence. Online detection is implemented with window length $\omega=10$. Figure \ref{solar_stat}(a) and Figure \ref{solar_stat}(b)
plot values of $\textsf{S}^3 \textsf{T}$ and the quadratic score statistic on a logarithmic scale, respectively. As we can observe, both statistics obtain peak detection statistics at around $t=227$, indicating both statistics can successfully detect the emergence of a solar flare.

\subsection{Water quality monitoring}

In this section, we consider a real-time water quality monitoring example for a sensor network deployed along a complex river system. The goal is to quickly detect contaminant spills that cause pollution of the river.

\begin{figure}[h!]
	\begin{center}
		\begin{tabular}{ccc}
			\includegraphics[width=.3\textwidth]{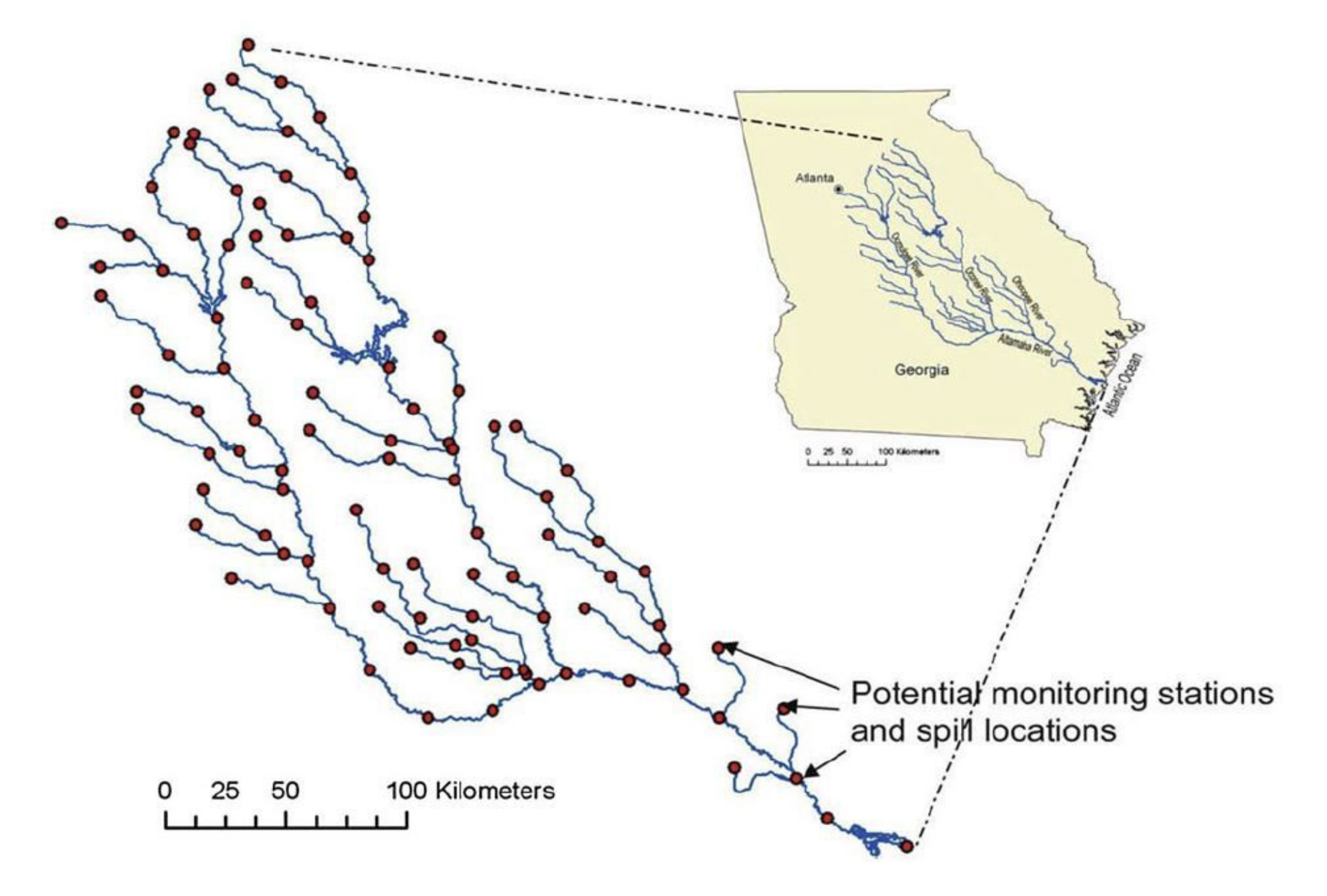} &
			\includegraphics[width=.3\textwidth]{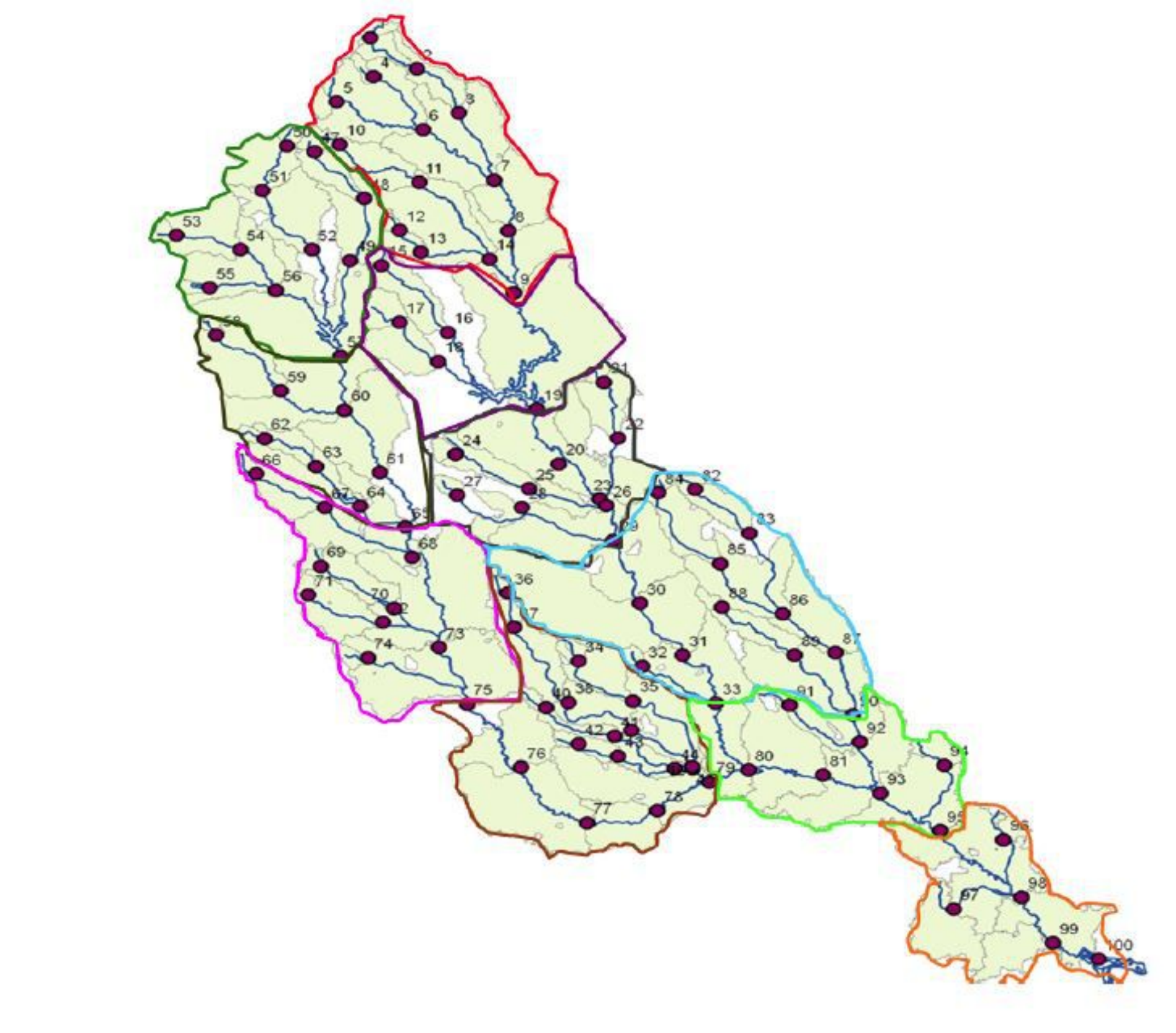} & 		\includegraphics[width=.4\textwidth]{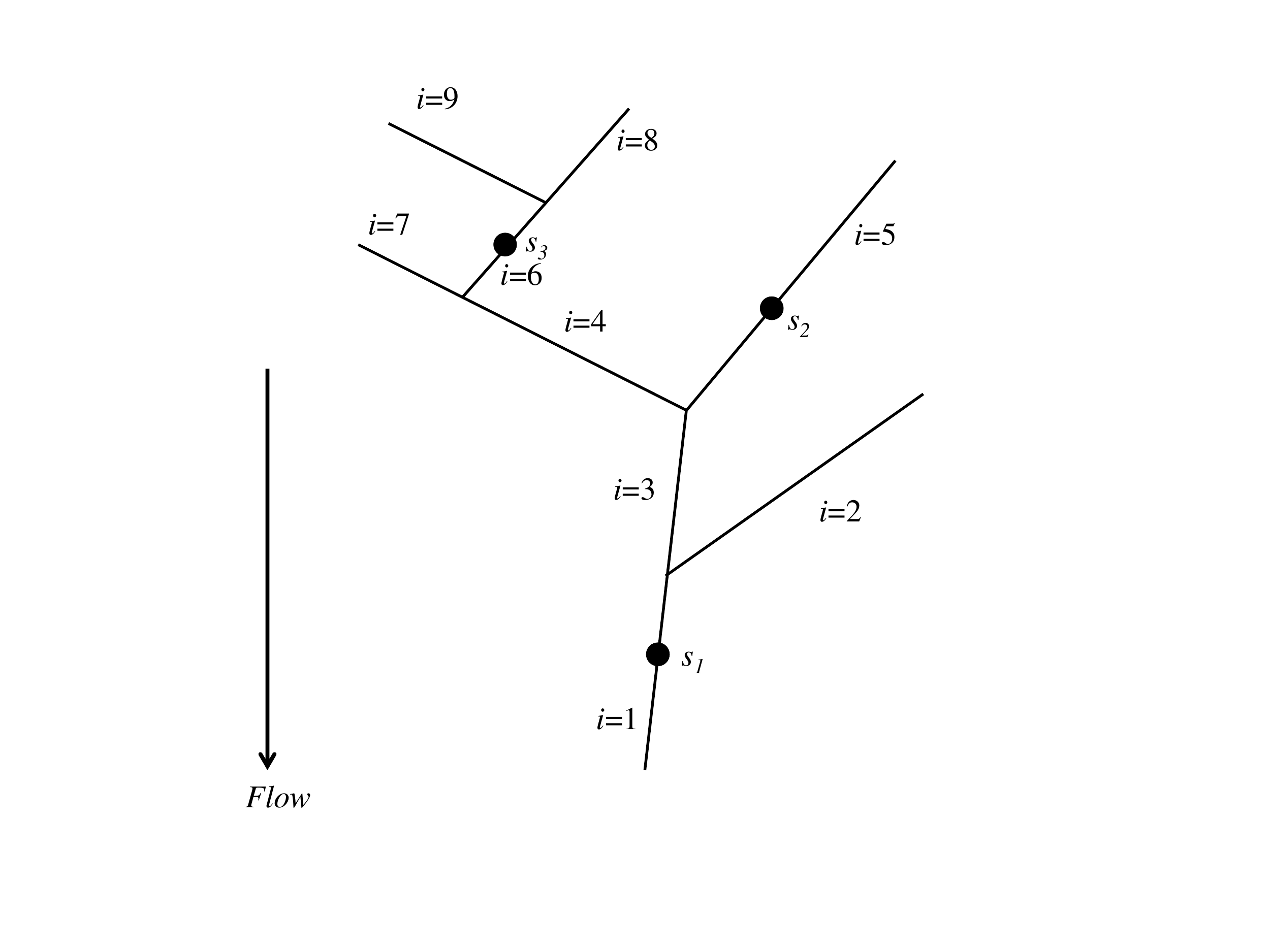}	
\\
			(a) & (b) & (c)
		\end{tabular}	
	\end{center}
	\caption{Water quality monitoring using sensor network: (a) shape and monitoring locations,  (b) ten sub-catchments of the Altamaha River \cite{telci2011contaminant}, and (c) an example of stream network with nine stream segments ($i=1,\cdots,9$) and three locations $s_1$, $s_2$, $s_3$.}
	\label{river_shape}
\end{figure}

We study the Altamaha River in Georgia, United States. The shape of the river is shown in Figure \ref{river_shape}(a). The nodes in the river network represent monitoring locations where concentration data is collected. 
The contaminant concentration data for such a river network is simulated by the Storm Water Management Model (SWMM) developed by the United States Environmental Protection Agency. SWMM requires geologic, geometric and fundamental hydrodynamics data to construct a river network. Given rainfall information, and the location, intensity, and duration of a contaminant spill, SWMM simulate the contaminant transport process through the river over a period. 
In the river dynamic simulation systems, rain events bring randomness to the contaminant transport. We use the same data as used in \citep{telci2011contaminant} to generate rain events. The Altamaha River watershed is divided into ten sub-catchments as shown in Figure \ref{river_shape}(b). The rainfall measurements are obtained from different United States Geological Survey stations close to these ten sub-catchments in 2006. Based on statistical analysis of these measurements, five rain patterns are generated for each sub-catchment. Each rain pattern describes time-dependent rainfall events and keeps changing hydrologic conditions in each-catchment during the simulation. Note that the rain patterns for each sub-catchment are different and thus there are $5^{10}$ possible combinations for the entire watershed.

Due to the nature of hydrodynamics, there is strong spatial correlation among the concentration data collected at different locations in the river network. However, the shape of the network and direction of the stream impose constraints on modeling such spatial correlation. For example, there should not be correlation for data collected at two locations that do not share a common flowing. A reasonable spatial correlation model is critical here.

We adopt the so-called ``tail-up'' spatial model for stream networks, which is proposed based on moving average constructions in \citep{ver2010moving}. The tail-up models have the following desired properties: (i) they use stream distance rather than the Euclidean distance, which is defined as the shortest distance along the stream network between two locations; (ii) statistical independence is imposed on the observations located on stream segments that do not share flowing water; (iii) proper weighting is incorporated on the entries of covariance matrix when the line segments in the network is splitting into multiple segments to ensure that the resulting covariance is stationary.

To explain the tail-up models, we first introduce some notations. A stream network consists of a finite number of stream segments and we index them with $i=1,2,\cdots$. Denote the index set of stream segment  as $I$, and the locations on the network as $s_j$, $j = 1,2,\cdots$. Let $D_{s_j} \subseteq I$ be the index set of all stream segments that are downstream of location $s_j$ (which means water from $s_j$ flows into these segments), including the segment containing $s_j$. Figure \ref{river_shape}(c) illustrates a simple stream network with $I = \{1,2,\cdots, 9\}$, $D_{s_1} = \{1\}$, $D_{s_2} = \{1,3,5\}$ and $D_{s_3} = \{1,3,4,6\}$.
Two locations, $s_j$ and $s_k$ are said to be ``flow-connected'' if $D_{s_j} \cap D_{s_k} = D_{s_j}$ or $D_{s_k}$. Finally, define
\begin{equation*}
B_{s_j, s_k} =  \left\{ \begin{array}{cl} \overline{ (D_{s_j} \cap D_{s_k}) } \cap (D_{s_j} \cup D_{s_k}), & \text{if $s_j$ and $s_k$ are flow-connected;} \\ \emptyset, & \text{otherwise.}
\end{array} \right.
\end{equation*}
Here $B_{s_j, s_k}$ is the set of stream segments between two locations, including the segment for the upstream location but excluding the one for the downstream location. For example, $B_{s_1, s_3} = \{3,4,6\}$ and $B_{s_2, s_3} = \emptyset$.
To ensure the stationarity of the variances, \cite{ver2010moving} suggests assigning weights to each stream segments in the network. In a stream network, one segment splits into two segments when it goes up-stream. For example, segment 1 splits into segments 2 and 3 in Figure \ref{river_shape}(c). One way to weight the segments is based on the flow volume of each segments. For example, we weight segments 2 and 3 by $w_2$ and $w_3$, where $w_2 + w_3 = 1$ and $w_2/w_3$ is equal to the ratio of the flow volume between segments 2 and 3.
Using tail-up models, the covariance between two locations, $s_j$ and $s_k$ on the stream network is given by
\begin{equation} \label{river_cov}
C(s_j, s_k | \bm{\zeta}) = \left\{ \begin{array}{ll}
0, & \text{if $s_j$ and $s_k$ are not flow-connected; } \\
\zeta_1, & \text{if $s_j = s_k$,}\\
\prod_{i \in B_{s_j, s_k}} \sqrt{w_i} \zeta_1 \rho\big(d(s_j, s_k)/\zeta_2 \big), & \text{otherwise;}
\end{array} \right.
\end{equation}
where $d(s_j, s_k)$ is the stream distance between $s_j$ and $s_k$,
$\zeta_1$ is the variance parameter, $\rho(\cdot | \zeta_2)$ is the correlation function with parameter $\zeta_2$, and $w_i$ is the weights on segment $i$. The correlation function $\rho(\cdot | \zeta_2)$ can be derived from many commonly used spatial models that we have discussed in Section \ref{formulation}.
For illustration, consider the example in Figure \ref{river_shape}(c). If an exponential model is used for spatial correlation, the covariance matrix of $s_1$, $s_2$ and $s_3$ can be constructed based on \eqref{river_cov} as follows,
\begin{equation*}
\left(
\begin{array}{ccc}
1 & \sqrt{w_3 w_5} & \sqrt{w_3 w_4 w_6} \\
\sqrt{w_3 w_5} & 1 & 0 \\
\sqrt{w_3 w_4 w_6} & 0 & 1
\end{array}
\right) \odot \left(
\begin{array}{ccc}
\zeta_0 + \zeta_1 & \zeta_1 e^{-d(s_1, s_2) /\zeta_2} & \zeta_1 e^{-d(s_1, s_3) /\zeta_2}  \\
\zeta_1 e^{-d(s_1, s_2) /\zeta_2}  & \zeta_0 + \zeta_1 & \zeta_1 e^{-d(s_2, s_3) /\zeta_2}  \\
\zeta_1 e^{-d(s_1, s_3) /\zeta_2}  & \zeta_1 e^{-d(s_2, s_3) /\zeta_2}  & \zeta_0 + \zeta_1
\end{array}
\right),
\end{equation*}
where $\odot$ denotes the Hadamard (element-wise) product operation between two matrices.

\begin{figure}[h!]
	\begin{center}
	\begin{tabular}{ccc}
		\includegraphics[width=.3\textwidth]{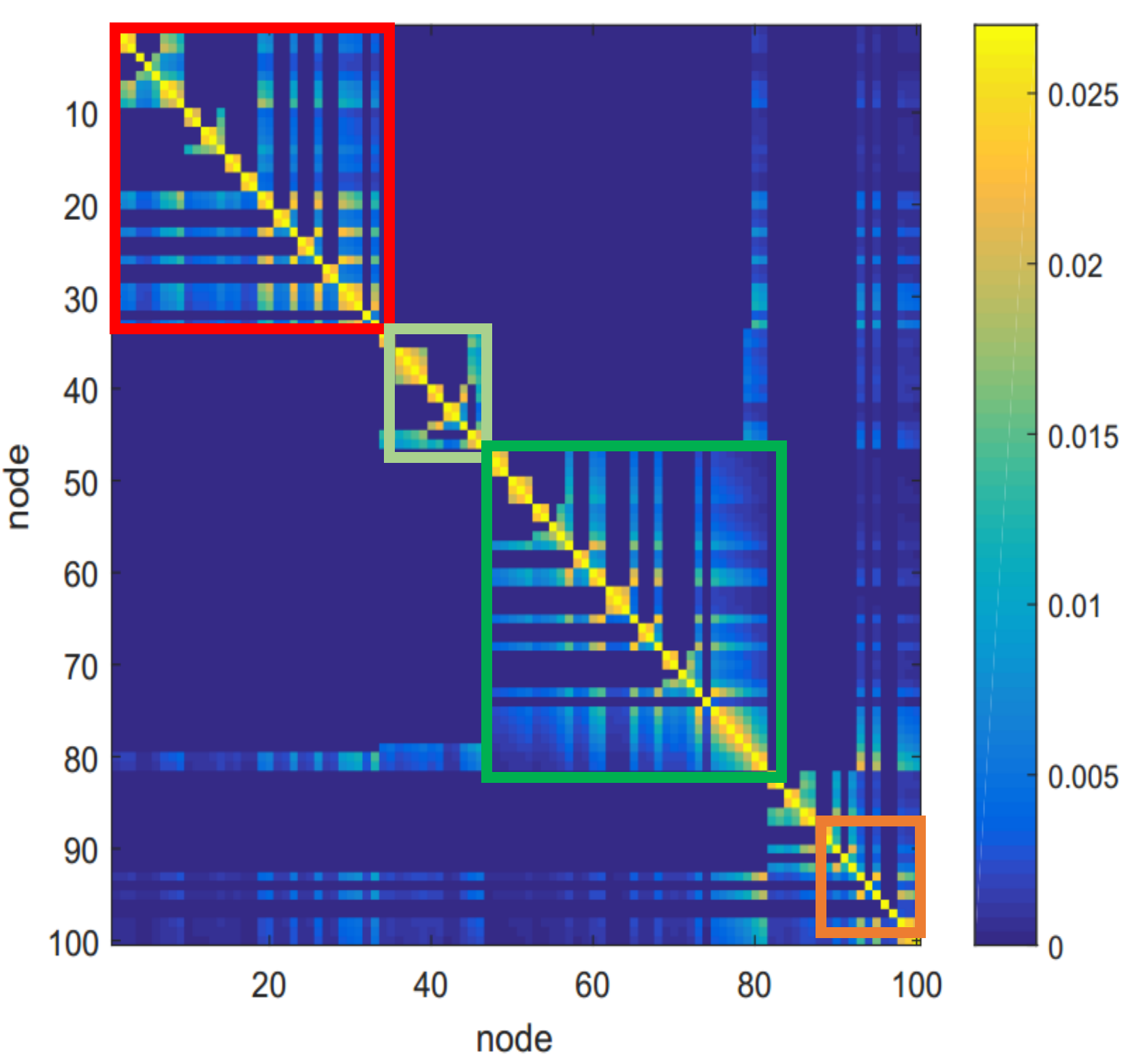} &
		\includegraphics[width=.25\textwidth]{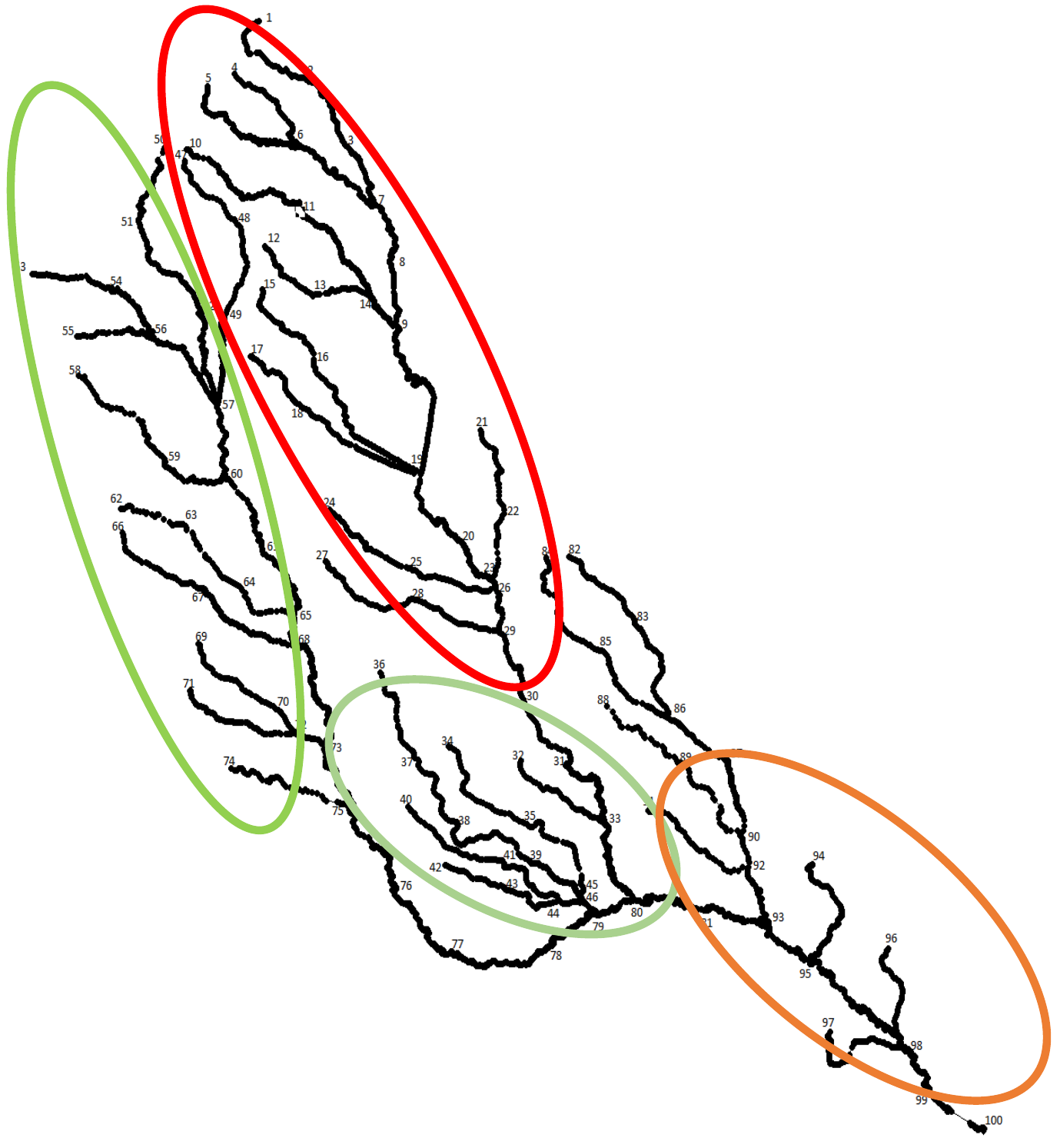}
		&		\includegraphics[width=.3\textwidth]{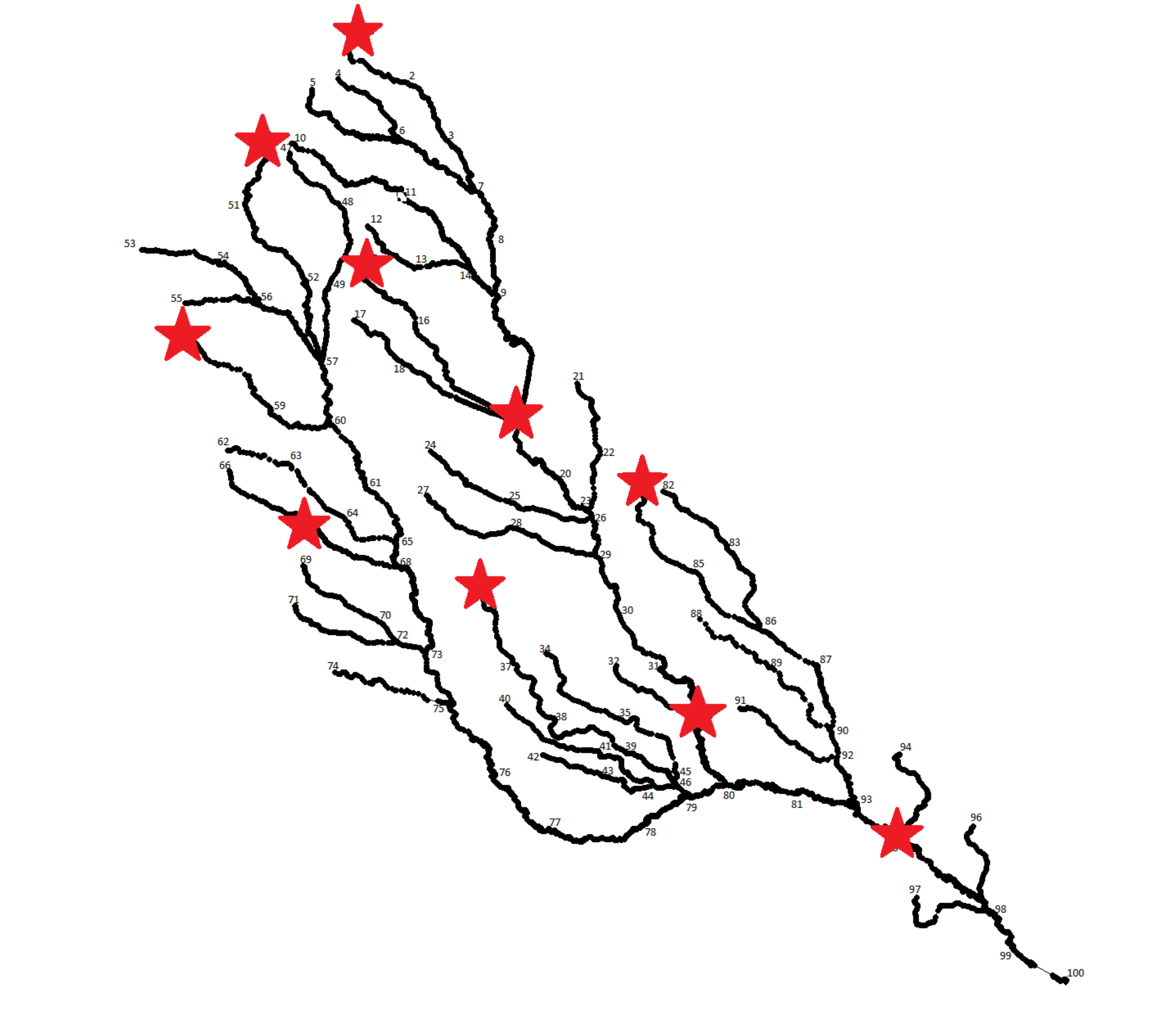}	
	\end{tabular}
	\end{center}
	\caption{(a) Visualization of the spatial covariance matrix using the tail-up model for 100 sensor network over the river system; the spatial covariance matrix has a block structure, with blocks in the matrix correspond to the branches of the river with matching colors in (b), and (c): nodes indexes of the Altamaha river network and  potential spill locations marked by red stars.}
	\label{Sw}
\end{figure}

In our case study, we use the tail-up model with exponential correlation function to model the data collected at different nodes on the Altamaha river network. The spatial covariance matrix for $p = 100$ nodes on the river network are constructed based on the stream distance and flow volume information. We use SWMM to generate in-control data and obtain the maximum likelihood estimator of the parameters in the model, $\hat{\zeta}_1 = 0.027$ and $\hat{\zeta}_2 = 0.68$. The covariance matrix is illustrated in Figure \ref{Sw}. 
For temporal correlation, we use a VAR(1) model $\bm{x}_{\ell} = \bm{\mu}_{x} + \theta \bm{x}_{\ell-1} + \bm{\epsilon}_{\ell}$,  where $\theta \in \mathbb{R}$ to capture the temporal correlation of a contaminant spill as suggested by \cite{Clement2007,Clement9}.

We apply the online change-point detection procedure based on $\textsf{S}^3 \textsf{T}$ to detect contaminant spills in the Altamaha river network. We also compare it with two other methods: (i) online detection based on the quadratic score statistic, and (ii) the Hotelling's $T^2$ chart.
Among the 100 nodes on the river network, 10 of them (nodes 1, 15, 19, 33, 36, 50, 58, 67, 84, 95, marked by red stars in Figure \ref{Sw}(c)) are used as possible contaminant spill locations, and the rest 90 nodes are used for collecting measurements every 15 minutes.
In each replication, we run SWMM to simulate the river network during a 10-day period. A single instantaneous spill with a spill location randomly selected from the ten possible locations is generated. The spill starting time is uniformly distributed between the first 15 to 20 hours. The intensity of the contaminant spills follows uniform distribution, and we consider three different levels: U(10, 100) (low), U(100, 250) (medium), and U(250, 500) (high) in units of gram/liter.
The thresholds for the three detection procedures are adjusted so that the in-control ARLs are 10 days (960 samples). For the two procedures based on $\textsf{S}^3 \textsf{T}$ and the quadratic score statistic, the length of the sliding window is chosen as 12.5 hours (50 samples). Table \ref{SWMM} reports the average and standard error of detection delays obtained from 100 simulated spills. For spills with high intensity, all three methods achieve similar performance in terms of detection delay, as strong signals are easier to be detected. However, when the signal is relatively weak (low and medium spill intensity), the proposed detection statistic $\textsf{S}^3 \textsf{T}$ significantly outperforms the other two methods.

\begin{table}[h!]
	\centering
	\caption{ Simulated expected detection delay in hours (Numbers in parentheses are standard errors).}
	\begin{tabular}{lccc}
		\toprule
		spill intensity &$\textsf{S}^3 \textsf{T}$    & Quadratic score statistic     & $T^2$ \\
		\midrule
		low   & \textbf{38.285 (3.655)} & 45.822 (4.675) & 52.959 (5.035) \\
		medium & \textbf{26.301 (1.679)} & 28.522 (1.873) & 30.753 (2.192) \\
		high  & 25.519 (1.697) & \textbf{25.489 (1.667)} & 25.563 (1.860) \\
		\bottomrule
	\end{tabular}%
	\label{SWMM}%
\end{table}%

\section{Conclusions}

In this paper, we propose a novel efficient score statistic $\textsf{S}^3 \textsf{T}$ to detect the emergence of a spatial-temporal signal from a noisy background. The statistic is able to jointly capture the spatial and temporal correlation and enjoys relatively low computational cost. An accurate approximation for its probability of a false alarm is presented. 
Numerical results on a simulated vector time series model and real applications show that the proposed statistic has a clear advantage when the signal is weak.

\section*{Acknowledgement}

The work is partially supported by NSF CMMI-1538746.



\bibliographystyle{plain}
\bibliography{ref}







\section*{Appendix A: Derivation of $\frac{\partial \ell}{\partial \bm{\mu}} \big|_{\bm{\mu}=\bm{0}, \gamma=0}$.}

The following propositions are used in the derivation of $\frac{\partial \ell}{\partial \bm{\mu}} \big|_{\bm{\mu}=\bm{0}, \gamma=0}$.

\begin{proposition} \label{prop1}
	Let $\bm{M}(t)$ be a nonsingular square matrix whose elements are functions of a scalar parameter $\alpha$. Then,
	$$
	\frac{\partial \bm{M}(t)^{-1}}{\partial \alpha} = - \bm{M}(t)^{-1} \frac{\partial \bm{M}(t)}{\partial \alpha} \bm{M}(t)^{-1}.
	$$
\end{proposition}

\begin{proposition} \label{prop2}
	Let $\bm{M}(t)$ be a nonsingular square matrix whose elements are functions of a scalar parameter $\alpha$. Then,
	$$
	\frac{\partial |\bm{M}(t)|}{\partial \alpha} = |\bm{M}(t)| \tr\bigg( \bm{M}(t)^{-1} \frac{\partial \bm{M}(t)}{\partial \alpha} \bigg).
	$$
\end{proposition}

By Proposition \ref{prop1}, we can calculate,
\begin{equation*}
\begin{aligned}
\frac{\log\big|\gamma \bm{V}_\tau(\theta) + \bm{\Sigma}_\tau\big|}{\partial \gamma} \Bigg| _{\bm{\mu}=\bm{0}, \gamma = 0} &= \frac{1}{\big|\gamma \bm{V}_\tau(\theta) + \bm{\Sigma}_\tau\big|} \big|\gamma \bm{V}_\tau(\theta) + \bm{\Sigma}_\tau\big| \tr\Big( (\gamma \bm{V}_\tau(\theta) + \bm{\Sigma}_\tau)^{-1} \bm{V}_\tau(\theta) \Big) \Bigg| _{\gamma = 0} \\
&= \tr\big(\bm{\Sigma}_\tau^{-1} \bm{V}_\tau(\theta)\big).
\end{aligned}
\end{equation*}

For convenience, here we use $\bm{y}$ and $\bm{\mu}$ to denote $\bm{y}_{(k+1:N)}$ and $\bm{\mu}_{(k+1:N)}$. By Proposition \ref{prop2}, we have,
\begin{equation*}
\begin{aligned}
\frac{\partial (\bm{y} - \bm{\mu})^\intercal(\gamma \bm{V}_\tau(\theta) + \bm{\Sigma}_\tau)^{-1} (\bm{y} - \bm{\mu})}{\partial \gamma} \Bigg| _{\bm{\mu}=\bm{0}, \gamma = 0} &= (\bm{y} - \bm{\mu})^\intercal \frac{\partial (\gamma \bm{V}_\tau(\theta) + \bm{\Sigma}_\tau)^{-1}}{\partial \gamma}  (\bm{y} - \bm{\mu}) \Bigg| _{\bm{\mu}=\bm{0}, \gamma = 0}\\
&= - \bm{y}^\intercal (\gamma \bm{V}_\tau(\theta) + \bm{\Sigma}_\tau)^{-1} \bm{V}_\tau(\theta) (\gamma \bm{V}_\tau(\theta) + \bm{\Sigma}_\tau)^{-1}  \bm{y} \Bigg| _{\gamma = 0} \\
&=- \bm{y}^\intercal \bm{\Sigma}_\tau^{-1} \bm{V}_\tau(\theta)\bm{\Sigma}_\tau^{-1}  \bm{y}.
\end{aligned}
\end{equation*}

Hence we have,
$$
\frac{\partial \ell}{\partial \bm{\mu}} \big|_{\bm{\mu}=\bm{0}, \gamma=0} = -\frac{1}{2}\tr\big(\bm{\Sigma}_\tau^{-1} \bm{V}_\tau(\theta)\big) + \frac{1}{2} \bm{y}^\intercal \bm{\Sigma}_\tau^{-1} \bm{V}_\tau(\theta) \bm{\Sigma}_\tau^{-1} \bm{y},
$$
as appeared in equation \eqref{score}.

\section*{Appendix B: Derivation of Equation \eqref{psi_fun}.}

Here we present the derivation of the cumulant generating function of $W(\tau, \theta)$ under the null hypothesis, i.e.\@ equation \eqref{psi_fun}.

Let $\bm{z} = \bm{\Sigma}_{\tau}^{-\frac{1}{2}} \bm{y}_{(k+1:N)}$. Under the null hypothesis, $\bm{z} \sim \mathcal{N}\big(\bm{0}, \bm{I}_{p\tau}\big)$. For convenience, here we use $\bm{B}$ to denote the $p\tau$ by $p\tau$ matrix $\bm{\Sigma}_\tau^{-\frac{1}{2}} \bm{V}_\tau(\theta) \bm{\Sigma}_\tau^{-\frac{1}{2}}$, and use $c$ and $d$ to denote $c(\tau, \theta)$ and $d(\tau, \theta)$, respectively. Then, we have
$$
W(\tau, \theta) = \frac{\bm{z}^\intercal \bm{B} \bm{z} - c}{\sqrt{d}}.
$$
Under the null hypothesis, the cumulant generating function of $W(\tau, \theta)$ can be calculated as
\begin{equation*}
\begin{aligned}
\psi(\xi) &= \log \E[\exp(\xi W(\tau, \theta))] = \log \E\Big[\exp\Big(\xi \Big( \frac{\bm{z}^\intercal \bm{B} \bm{z} - c}{\sqrt{d}}\Big)\Big)\Big]\\
& = - \xi \frac{c}{\sqrt{d}} + \log \E\Big[\exp\Big(\frac{\xi \bm{z}^\intercal \bm{B} \bm{z}}{\sqrt{d}}\Big)\Big] \\
& = - \xi \frac{c}{\sqrt{d}} + \log \int_{\bm{z}} \exp\Big(\frac{\xi \bm{z}^\intercal \bm{B} \bm{z}}{\sqrt{d}}\Big) \frac{1}{(2\pi)^{\frac{p\tau}{2}}} \exp\Big(-\frac{1}{2}\bm{z}^\intercal\bm{z}\Big) d\bm{z} \\
& = - \xi \frac{c}{\sqrt{d}} + \log \int_{\bm{z}} \frac{1}{(2\pi)^{\frac{p\tau}{2}}} \exp\Big(-\frac{1}{2}\bm{z}^\intercal \Big(\bm{I}_{p\tau} - \frac{2\xi \bm{B}}{\sqrt{d}}\Big) \bm{z}\Big) d\bm{z} \\
& = - \xi \frac{c}{\sqrt{d}} + \log \Big|\bm{I}_{p\tau} - \frac{2\xi \bm{B}}{\sqrt{d}} \Big|^{-\frac{1}{2}},
\end{aligned}
\end{equation*}
which is equivalent to equation \eqref{psi_fun}. Note that the last equation uses the fact that
$$
\int_{\bm{z}} \frac{1}{(2\pi)^{\frac{p\tau}{2}}} \exp\Big(-\frac{1}{2}\bm{z}^\intercal \Big(\bm{I}_{p\tau} - \frac{2\xi \bm{B}}{\sqrt{d}}\Big) \bm{z}\Big) d\bm{z} = \Big|\bm{I}_{p\tau} - \frac{2\xi \bm{B}}{\sqrt{d}} \Big|^{-\frac{1}{2}}.
$$

\section*{Appendix C: Proof of Theorem \ref{theorem1}}

After discretizing the parameter space, $W(\tau, \theta)$ can be treated as a two-dimensional Gaussian random field, which can be completely characterized by its covariance function. The following lemma computes the covariance function of $W(\tau, \theta)$.
\begin{lemma} \label{cov1}
	Under the null hypothesis, the covariance function of $W(\tau, \theta)$ is
	\begin{equation} \label{cov function}
	\begin{aligned}
	\Cov[& W(n, \theta_{1}),  W(m, \theta_{2})] = \frac{\tr\big(\bm{A}_{n}(\theta_{1}) \bm{A}_{n}(\theta_{2})\big)}{\Big[ \tr\big( \bm{A}_{n}(\theta_{1})\bm{A}_{n}(\theta_{1} ) \big) \tr\big(  \bm{A}_{m}(\theta_{2})\bm{A}_{m}(\theta_{2} ) \big) \Big]^{1/2}},
	\end{aligned}
	\end{equation}
	where $n \leq m$.
\end{lemma}

The following lemma shows that the first order approximation of the covariance function in \eqref{cov function} does not have any cross product term. Thus, the two-dimensional random field is further decomposed as a sum of two independent one-dimensional random processes.
\begin{lemma} \label{taylor1}
	Assuming that $\delta$ and $i \in Z$ are small relative to $\theta$ and $\tau$, respectively, the first order approximation of the covariance function in \eqref{cov function} is given as,
	\begin{equation} \label{taylor1_main}
	\begin{aligned}
	&\Cov[ W(\tau, \theta), W(\tau + i, \theta + \delta) ]  \approx 1 - \gamma^{2}(\tau, \theta) \delta^{2} - \frac{\mu(\tau, \theta)}{2\tau}i + o(\delta^{2}) + o(i),
	\end{aligned}
	\end{equation}
	where
	\begin{equation} \label{gamma_fun}
	\gamma(\tau, \theta) = \frac{\tr\big(\dot{\bm{A}}_{\tau}(\theta)\bm{A}_{\tau}(\theta)\big)}{\tr\big(\bm{A}_{\tau}(\theta)\bm{A}_{\tau}(\theta)\big)},
	\end{equation}
	$\mu(\tau, \theta)$ is defined in \eqref{mu_fun}, and $\dot{\bm{A}}_{\tau}(\theta) = \partial \bm{A}_{\tau}(\theta)/\partial \theta$.
\end{lemma}

The following two Lemmas are needed in the proof. Both Lemmas are proved in \cite{xie2012spectrum}.

\begin{lemma} \label{decomposition}
	Assume $\xi \rightarrow \infty$, $b \rightarrow \infty$, $N \rightarrow \infty$, with $\frac{\xi}{b} \approx 1$ and $\frac{b}{N} \approx c$, where $c > 0 $ is some constant. The discretized process $b \Big[W\Big( \tau + i, \theta + \frac{\Delta}{\sqrt{N}j} \Big) - \xi\Big]$, where $i$ is an integer and $j \geq 0$, conditioned on $W(\tau, \theta) = \xi$ can be written as a sum of two independent processes:
	\[
	\bigg\{ b \Big[W\Big( \tau + i, \theta + \frac{\Delta}{\sqrt{N}}j \Big) - \xi\Big] \bigg| W(\tau, \theta) = \xi \bigg\} = S_{i} + V_{j},
	\]
	where $S_{i} = \sum_{l = 1}^{i} a_{\ell}$ with
	\begin{equation*}
	a_{\ell} \sim N\bigg( -\frac{\mu(\tau, \theta)}{2\tau}b^{2}, \frac{\mu(\tau, \theta)}{\tau}b^{2} \bigg),
	\end{equation*}
	and
	\begin{equation*}
	V_{j} = \sqrt{2} \gamma(\tau, \theta) \frac{b}{\sqrt{N}}\Delta j V - \gamma^{2}(\tau, \theta)\frac{b^{2}}{N} \Delta^{2}j^{2},
	\end{equation*}
	with $V\sim N(0, 1)$. $\mu(\tau, \theta)$ and $\gamma(\tau, \theta)$ are defined in \eqref{mu_fun} and \eqref{gamma_fun}, respectively.
	
\end{lemma}

\begin{lemma} \label{shrinking}
	Assume $x_{1},\ x_{2}, \cdots$ are i.i.d.\@ $N(-\mu_{1}, \sigma_{1}^{2})$ random variables ($\mu_{1} > 0$). Define the random walk $S_{0} = 0$, $S_{i} = \sum_{l=1}^{i} x_{l}, \ i=1,2,\cdots$, and the smooth varying random process $V_{j} = \beta \Delta j V - \frac{\beta^{2}}{2} \Delta^{2} j^{2}$, for some constants $\Delta > 0$, $\beta > 0$. As $\Delta \rightarrow 0$, for some constant $\alpha$, we have
	\begin{equation*}
	\begin{aligned}
	&\frac{1}{\Delta} \int_{0}^{\infty} e^{-\alpha x} \mathbb P\bigg( \max_{i \geq 1} S_{i} \leq -x \bigg) \mathbb P\bigg( \max_{i \leq 0} S_{i} + \max_{j \geq 1} V_{j} \leq -x \bigg)dx \xrightarrow[]{\Delta \rightarrow 0} \frac{|\beta|}{\sqrt{2\pi}} \bigg(\frac{2\mu_{1}^{2}}{\sigma_{1}^{2}} \bigg) \nu\bigg( \frac{2\mu_{1}}{\sigma_{1}} \bigg),
	\end{aligned}
	\end{equation*}
	where $\nu(x)$ is defined in \eqref{nu_fun}.
	
\end{lemma}

In the following, we go through the main steps that lead to the approximation of the false alarm rate in Theorem \ref{theorem1} for the case of $d = 1$.

\noindent{\bf Step 1:} We first discretize the parameter $\theta \in [\theta_{1}, \theta_{2}]$ by a rectangular mesh grid of size $\frac{\Delta}{\sqrt{N}}$, where $\Delta > 0$ is a small number. Note that the discretization mentioned here is used for asymptotic analysis only. The probability of false alarm can be approximated as
\begin{equation} \label{Pe}
\mathbb P\Big( \max_{(i,j) \in D} W\Big(i, j\frac{\Delta}{\sqrt{N}} \Big) \geq b \Big),
\end{equation}
where $D$ is the index set
\begin{equation*}
D = \Big\{ (i,j): 1 \leq i \leq N,\ \theta_{1} \leq j\frac{\Delta}{\sqrt{N}} \leq \theta_{2} \Big\},
\end{equation*}
which covers the entire parameter space. Let $J(i_{0}, j_{0})$ denote everything to the ``future'' of the current index $(i_{0}, j_{0})$ in the parameter space, i.e.\@,
\begin{equation*}
J(i_{0}, j_{0}) = \{ (i, j) \in D: j \geq j_{0}, \text{or}\ i \geq i_{0}\ \text{and}\ j = j_0 \}.
\end{equation*}

Using the similar approach as in (\cite{siegmund1988approximate}), the event $\Big\{\max_{(i,j) \in D} W\Big(i, j\frac{\Delta}{\sqrt{N}} \Big) \geq b \Big\}$ can be decomposed into a series of ``last hitting events'', for which $(i_0, j_0)$ is the ``last'' location where $ W\Big(i, j\frac{\Delta}{\sqrt{N}} \Big)$ hits the threshold $b$. Then, the probability in \eqref{Pe} can be written as the sum of probabilities of $ W\Big(i, j\frac{\Delta}{\sqrt{N}} \Big)$ last hits $b$ at $(i_0, j_0)$ over all possible $(i_0, j_0)$:
\begin{equation} \label{last_hit}
\begin{aligned}
&\mathbb P\Big( \max_{(i,j) \in D} W\Big(i, j\frac{\Delta}{\sqrt{N}} \Big) \geq b \Big) \approx \sum_{(i_0,j_0) \in D} \mathbb P \Big(  W\Big(i_{0}, j_{0}\frac{\Delta}{\sqrt{N}} \Big) \geq b, \max_{(i,j) \in J(i_{0}, j_{0})} W\Big(i, j\frac{\Delta}{\sqrt{N}} \Big) < b \Big) \\
& = \sum_{(i_0,j_0) \in D} \int_{0}^{\infty} \mathbb P \Big(  W\Big(i_{0}, j_{0}\frac{\Delta}{\sqrt{N}} \Big) = b + \frac{x}{b} \Big)\\
& \quad \cdot \mathbb P\Big(\max_{(i,j) \in J(i_{0}, j_{0})} W\Big(i, j\frac{\Delta}{\sqrt{N}} \Big) < b \Big| W\Big(i_{0}, j_{0}\frac{\Delta}{\sqrt{N}} \Big) = b + \frac{x}{b} \Big)\frac{dx}{b}.
\end{aligned}
\end{equation}

\noindent{\bf Step 2:} In the following, we obtain an approximation on the probability $ \mathbb P \Big(  W\Big(i_{0}, j_{0}\frac{\Delta}{\sqrt{N}} \Big) = b + \frac{x}{b} \Big)\frac{dx}{b} $. To simplify the notation, we denote $W\Big(i_{0}, j_{0}\frac{\Delta}{\sqrt{N}} \Big)$ as $W$ here. The key idea is to approximate $W$ as a Gaussian random field. The Gaussian approximation performs well when the probability of interest is close to the mean of the true distribution, but suffers from deviation if the probability is in the tail of the true distribution. Hence, we apply the change-of-measure technique to shift the mean of the random field $W$ to the threshold $b$.

Denote the cumulant generating function of $W$ as $\psi(\xi) = \log \E[\exp(\xi W)]$. To construct the new probability measure, we first choose a $\xi_{0} > 0$ such that $\psi'(\xi) = b$. The new probability measure $dF_{\xi_{0}}$ is constructed using exponential embedding, as follows
\begin{equation*}
dF_{\xi_{0}} = \exp \big( \xi_{0}W - \psi(\xi_{0}) \big)dF,
\end{equation*}
where $dF$ is the original distribution of $W$. Let $\E_{\xi_{0}}$ and $\mathbb P_{\xi_{0}}$ denote the expectation and probability under the new measure $dF_{\xi_{0}}$, respectively. It can be verified that under the new measure
\begin{equation*}
\begin{aligned}
\E_{\xi_{0}}[W] &= \E\big[ W \exp \big( \xi_{0}W - \psi(\xi_{0}) \big) \big]  = e^{-\psi(\xi_{0})} \frac{\partial e^{\psi(\xi)}}{\partial \xi} \bigg|_{\xi = \xi_{0}} = \psi'(\xi) = b,
\end{aligned}
\end{equation*}
namely, the mean of $W$ is close to the threshold $b$ under the new probability measure.

The threshold crossing probability can be rewritten as
\begin{equation} \label{change-of-measure}
\begin{aligned}
\mathbb P \Big( W = b + \frac{x}{b} \Big) &= \E_{\xi_{0}} \bigg[ \frac{1}{\exp [\xi_{0} W - \psi(\xi_{0})]}\mathbbm{1}\{\ W = b + \frac{x}{b} \}\bigg] \\
& = \exp\Big[ \psi(\xi_{0}) - \xi_{0}\Big(b + \frac{x}{b}\Big) \Big] \mathbb P_{\xi_{0}} \Big(W = b + \frac{x}{b}\Big).
\end{aligned}
\end{equation}

Now we can apply the Gaussian approximation to obtain $\mathbb P_{\xi_{0}} \Big(W = b + \frac{x}{b}\Big)$ and use \eqref{change-of-measure} to get the original probability. By treating $W$ as a normal random variable with mean $b$ and variance $\sigma^{2}_{\xi_{0}}$, we have
\begin{equation*}
\begin{aligned}
\mathbb P_{\xi_{0}} \Big(W = b + \frac{x}{b}\Big) &= \frac{1}{\sqrt{2\pi}\sigma_{\xi_{0}}} \exp \bigg( \frac{-x^{2}}{2b^{2} \sigma^{2}_{\xi_{0}}} \bigg) \approx \frac{1}{\sqrt{2\pi}\sigma_{\xi_{0}}}.
\end{aligned}
\end{equation*}
Note that in \eqref{last_hit}, the integrands with smaller $x$ values contribute more to the integration, since the integrand decays exponentially fast with $x$. Now, when $b \rightarrow \infty$, $\frac{x}{b} \rightarrow 0$ for small $x$, and hence $\exp \bigg( \frac{-x^{2}}{2b^{2} \sigma^{2}_{\xi_{0}}} \bigg) \rightarrow 1$. The above argument is similar to those used for Laplace's method.

The cumulant generating function of $W$ can be calculated as
\begin{equation*}
\begin{aligned}
\psi(\xi) &= -\xi \frac{\tr\big( \bm{\Sigma}_{\tau}^{-1} \bm{V}_{\tau}(\theta) \big)}{ \Big[2 \tr\big(\bm{\Sigma}_{\tau}^{-1} \bm{V}_{\tau}(\theta) \bm{\Sigma}_{\tau}^{-1} \bm{V}_{\tau}(\theta)\big)\Big]^{1/2}} - \frac{1}{2} \log \bigg| \bm{I}_{p\tau} - \frac{2\xi \bm{\Sigma}_{\tau}^{1/2} \bm{V}_{\tau}(\theta)\bm{\Sigma}_{\tau}^{1/2}}{ \Big[2 \tr\big(\bm{\Sigma}_{\tau}^{-1} \bm{V}_{\tau}(\theta) \bm{\Sigma}_{\tau}^{-1} \bm{V}_{\tau}(\theta)\big)\Big]^{1/2}} \bigg|.
\end{aligned}
\end{equation*}
Hence $\xi_{0}$ can be obtained by solving the following equation numerically,

\begin{equation*}
\begin{aligned}
\frac{1}{\sqrt{d(\tau, \theta)}}\tr\Big( \Big[\bm{I}_{p\tau} - \frac{2\xi_{0} \bm{B}_{\tau}(\theta) }{ \sqrt{d(\tau, \theta)} }\Big]^{-1} \bm{B}_{\tau}(\theta) - \bm{A}_{\tau}(\theta) \Big) = b.
\end{aligned}
\end{equation*}

Eventually, we have
\begin{equation} \label{first_part}
\mathbb P \Big(  W\Big(i_{0}, j_{0}\frac{\Delta}{\sqrt{N}} \Big)  = b + \frac{x}{b} \Big) \approx g\Big(i_{0}, j_{0}\Big) \exp\Big(-\frac{\xi_0}{b}x\Big),
\end{equation}
where $g()$ follows the definition in \eqref{g_fun}.

\noindent{\bf Step 3:} Next we tackle with the conditional probability $ \mathbb P\Big(\max_{(i,j) \in J(i_{0}, j_{0})} W\Big(i, j\frac{\Delta}{\sqrt{N}} \Big) < b \Big| W\Big(i_{0}, j_{0}\frac{\Delta}{\sqrt{N}} \Big) = b + \frac{x}{b} \Big)$.
The first order expansion of the covariance function given by Lemma \ref{taylor1} does not have any cross product term, which implies that if we approximate $W(\tau, \theta)$ as a Gaussian random field, it can be decomposed as a sum of two independent one dimensional random processes.

By Lemma \ref{decomposition}, the conditional probability can be written in terms of the decomposed random processes using the techniques in \cite{siegmund1988approximate} and \cite{kim1989likelihood} as follows,
\begin{equation} \label{conditional}
\begin{aligned}
\mathbb P\Big(&\max_{(i,j) \in J(i_{0}, j_{0})} W\Big(i, j\frac{\Delta}{\sqrt{N}} \Big) < b \Big| W\Big(i_{0}, j_{0}\frac{\Delta}{\sqrt{N}} \Big) = b + \frac{x}{b} \Big) \\
&= \mathbb P \bigg( \max_{(i,j) \in J(i_{0}, j_{0})} b\bigg[ W\Big(i, j\frac{\Delta}{\sqrt{N}} \Big) -W\Big(i_{0}, j_{0}\frac{\Delta}{\sqrt{N}} \Big)  \bigg] \leq -x \bigg| W\Big(i_{0}, j_{0}\frac{\Delta}{\sqrt{N}} \Big) = b + \frac{x}{b}  \bigg) \\
&\approx \mathbb P\bigg( \max_{i \geq 1} S_{i} \leq -x \bigg) \mathbb P \bigg( \max_{i \leq 0} S_{i} + \max_{j \geq 1} V_{j} \leq -x \bigg).
\end{aligned}
\end{equation}

(4) Combine the approximations in \eqref{first_part} and \eqref{conditional}, the approximated false alarm rate becomes,
\begin{equation} \label{approx}
\begin{aligned}
&\mathbb P\Big( \max_{(i,j) \in D} W\Big(i, j\frac{\Delta}{\sqrt{N}} \Big) \geq b \Big) \\ &\approx \sum_{(i_0,j_0) \in D} g\Big( i_{0}, j_{0} \frac{\Delta}{\sqrt{N}} \Big) \frac{\Delta}{\sqrt{N}} \frac{\sqrt{N}}{\Delta b} \int_{0}^{\infty} \exp\Big( -\frac{\xi_{0}}{b} x \Big)\\
&\quad \cdot\mathbb P\bigg( \max_{i \geq 1} S_{i} \leq -x \bigg) \mathbb P \bigg( \max_{i \leq 0} S_{i} + \max_{j \geq 1} V_{j} \leq -x \bigg) dx.
\end{aligned}
\end{equation}

Lemma \ref{shrinking} enables us to find an expression for the integration in \eqref{approx}.

Finally, by Lemma \ref{shrinking} with $\alpha = \frac{\xi_{0}}{b}$, $\beta = \sqrt{2} \gamma(\tau, \theta) \frac{b}{\sqrt{N}}$, $\mu_{1} = \frac{\mu(\tau, \theta)}{2\tau} b^{2}$ and $\sigma_{1}^{2}  = \frac{\mu(\tau, \theta)}{\tau} b^{2}$, we have the approximated significance level
\begin{equation} \label{Riemann_sum}
\begin{aligned}
&\frac{1}{2\sqrt{\pi}} \sum_{(i_{0}, j_{0})\in D} g \bigg( i_{0}, j_{0}\frac{\Delta}{\sqrt{N}} \bigg) \frac{b^{2}\mu(i_{0}, j_{0}\frac{\Delta}{\sqrt{N}})}{N -i_{0}} \cdot \nu\Bigg( \sqrt{\frac{b^{2} \mu(i_{0}, j_{0} \frac{\Delta}{\sqrt{N}})}{N - i_{0}}} \Bigg) \gamma\bigg( i_{0}, j_{0} \frac{\Delta}{\sqrt{N}} \bigg)\frac{\Delta}{\sqrt{N}}.
\end{aligned}
\end{equation}
As $\Delta \rightarrow 0$, the Riemann sum \eqref{Riemann_sum} converges to the approximation in Theorem \ref{theorem1}.

\section*{Appendix D: Proof of Lemma \ref{cov1}.}

\begin{proof}
	
	
	Let $\bm{C}_{\tau}(\theta) = \bm{\Sigma}_{\tau}^{-1} \bm{V}_{\tau}(\theta)\bm{\Sigma}_{\tau}^{-1}$, and rewrite $\bm{C}_{m}(\theta_{2})$ as,
	$$
	\bm{C}_{m}(\theta_{2}) = \begin{bmatrix}
	\bm{C}_{11}(\theta_{2}) & \bm{C}_{12}(\theta_{2}) \\
	\bm{C}_{21}(\theta_{2}) & \bm{C}_{n}(\theta_{2})
	\end{bmatrix}.
	$$
	Denote $\bm{y}_{(T-\tau+1:T)}$ as $Y_{\tau}$, and let
	$$
	Y_{m} = \begin{bmatrix}
	Y_{\Delta} \\ Y_{n}
	\end{bmatrix}.
	$$
	We have,
	\begin{equation} \label{lamma11}
	\Cov[ W(n, \theta_{1}),  W(m, \theta_{2})] =  \frac{ \E[Y_{n}^\intercal\bm{C}_{n}(\theta_{1})Y_{n} Y_{m}^\intercal\bm{C}_{m}(\theta_{2})Y_{m} ] - \E[Y_{n}^\intercal\bm{C}_{n}(\theta_{1})Y_{n}] \E[Y_{m}^\intercal\bm{C}_{m}(\theta_{2})Y_{m}]  }{2( \tr\{ \bm{A}_{n}(\theta_{1})\bm{A}_{n}(\theta_{1} ) \} \tr\{  \bm{A}_{m}(\theta_{2})\bm{A}_{m}(\theta_{2} ) \} )^{1/2}}.
	\end{equation}
	
	The first term in the numerator is,
	\begin{equation} \label{lamma12}
	\begin{aligned}
	&\E[Y_{n}^\intercal\bm{C}_{n}(\theta_{1})Y_{n} Y_{m}^\intercal\bm{C}_{m}(\theta_{2})Y_{m} ] \\ &= \E[ (Y_{n}^\intercal\bm{C}_{n}(\theta_{1})Y_{n}) ( Y_{\Delta}^\intercal\bm{C}_{11}(\theta_{2})Y_{\Delta} + Y_{n}^\intercal\bm{C}_{n}(\theta_{2})Y_{n} + Y_{n}^\intercal\bm{C}_{21}(\theta_{2})Y_{\Delta} + Y_{\Delta}^\intercal\bm{C}_{12}(\theta_{2})Y_{n} ) ] \\
	& = \E[Y_{n}^\intercal\bm{C}_{n}(\theta_{1})Y_{n} Y_{n}^\intercal\bm{C}_{n}(\theta_{2})Y_{n} ] + \E[Y_{n}^\intercal\bm{C}_{n}(\theta_{1})Y_{n}] \E[Y_{\Delta}^\intercal\bm{C}_{11}(\theta_{1})Y_{\Delta}] \\
	& = 2\tr\{ \bm{A}_{n}(\theta_{1})\bm{A}_{n}(\theta_{2}) \} +  \tr\{\bm{A}_{n}(\theta_{1})\} \tr\{\bm{A}_{n}(\theta_{2})\} + \E[Y_{n}^\intercal\bm{C}_{n}(\theta_{1})Y_{n}] \E[Y_{\Delta}^\intercal\bm{C}_{11}(\theta_{1})Y_{\Delta}].
	\end{aligned}
	\end{equation}
	Note that we have utilized the fact that under null hypothesis, $Y_{\Delta}$ and $Y_{n}$ are independent and $\E[Y_{\Delta}] = 0$.
	
	The second term in the numerator is,
	\begin{equation} \label{lamma13}
	\begin{aligned}
	&\E[Y_{n}^\intercal\bm{C}_{n}(\theta_{1})Y_{n}] \E[Y_{m}^\intercal\bm{C}_{m}(\theta_{2})Y_{m}] \\ &= \E[ Y_{n}^\intercal\bm{C}_{n}(\theta_{1})Y_{n}] \E[  Y_{\Delta}^\intercal\bm{C}_{11}(\theta_{2})Y_{\Delta} + Y_{n}^\intercal\bm{C}_{n}(\theta_{2})Y_{n} + Y_{n}^\intercal\bm{C}_{21}(\theta_{2})Y_{\Delta} + Y_{\Delta}^\intercal\bm{C}_{12}(\theta_{2})Y_{n}  ] \\
	& = \E[Y_{n}^\intercal\bm{C}_{n}(\theta_{1})Y_{n}] \E[ Y_{n}^\intercal\bm{C}_{n}(\theta_{2})Y_{n} ] + \E[Y_{n}^\intercal\bm{C}_{n}(\theta_{1})Y_{n}] \E[Y_{\Delta}^\intercal\bm{C}_{11}(\theta_{1})Y_{\Delta}] \\
	& = \tr\{\bm{A}_{n}(\theta_{1})\} \tr\{\bm{A}_{n}(\theta_{2})\} + \E[Y_{n}^\intercal\bm{C}_{n}(\theta_{1})Y_{n}] \E[Y_{\Delta}^\intercal\bm{C}_{11}(\theta_{1})Y_{\Delta}].
	\end{aligned}
	\end{equation}
	
	By combining \eqref{lamma11}, \eqref{lamma12}, \eqref{lamma13}, we obtain the covariance function in Lamma \ref{cov1}.
	
\end{proof}

\section*{Appendix E: Proof of Lemma \ref{taylor1}.}
\begin{proof}
	We approximate the covariance function by expanding each term in \eqref{cov function} at $\theta$ and keeping only the first order terms.
	
	The numerator in \eqref{cov function} is approximated as,
	\begin{equation} \label{lemma2_num}
	\begin{aligned}
	\tr \big( \bm{A}_{\tau}(\theta + \delta) \bm{A}_{\tau}(\theta) \big) &\approx \tr \big( \bm{A}_{\tau}(\theta) \bm{A}_{\tau}(\theta) \big) + \delta \tr \big( \dot{\bm{A}}_{\tau}(\theta) \bm{A}_{\tau}(\theta) \big) \\
	&= \tr \big( \bm{A}_{\tau}(\theta) \bm{A}_{\tau}(\theta) \big) (1 + \delta \gamma(\tau, \theta)).
	\end{aligned}
	\end{equation}
	
	Partition the matrix $\bm{A}_{\tau + i} (\theta + \delta)$ as the follows,
	\begin{equation*}
	\bm{A}_{\tau + i} (\theta + \delta) = \begin{bmatrix}
	\bm{A}_{11}(\theta + \delta) & \bm{A}_{12}(\theta + \delta) \\
	\bm{A}_{21}(\theta + \delta) & \bm{A}_{\tau}(\theta + \delta)
	\end{bmatrix}.
	\end{equation*}
	Then rewrite the second term in the denominator in \eqref{cov function} as,
	\begin{equation*}
	\begin{aligned}
	&\tr\big( \bm{A}_{\tau + i}(\theta + \delta)\bm{A}_{\tau + i}(\theta + \delta) \big) \\
	&= \tr\big( \bm{A}_{11}(\theta + \delta)\bm{A}_{11}(\theta + \delta) \big) + \tr\big( \bm{A}_{12}(\theta + \delta)\bm{A}_{21}(\theta + \delta) \big) \\ &+\tr\big( \bm{A}_{21}(\theta + \delta)\bm{A}_{12}(\theta + \delta) \big) +\tr\big( \bm{A}_{\tau}(\theta + \delta)\bm{A}_{\tau}(\theta + \delta) \big).
	\end{aligned}
	\end{equation*}
	
	After expanding each term at $\theta$, the denominator in \eqref{cov function} can be approximated as,
	\begin{equation} \label{lemma2_den}
	\Big[ \tr\big( \bm{A}_{\tau}(\theta)\bm{A}_{\tau}(\theta) \big) \tr\big(  \bm{A}_{\tau + i}(\theta + \delta)\bm{A}_{\tau + i}(\theta + \delta ) \big) \Big]^{1/2} \approx \tr \big( \bm{A}_{\tau}(\theta) \bm{A}_{\tau}(\theta) \big) \sqrt{1 + 2\delta a} \sqrt{1 + b},
	\end{equation}
	where
	\begin{equation} \label{lemma2_a}
	a = \frac{\tr\big( \dot{\bm{A}}_{11}(\theta)\bm{A}_{11}(\theta) \big) + \tr\big( \dot{\bm{A}}_{12}(\theta)\bm{A}_{21}(\theta) \big) +\tr\big( \dot{\bm{A}}_{21}(\theta)\bm{A}_{12}(\theta) \big) +\tr\big( \dot{\bm{A}}_{\tau}(\theta)\bm{A}_{\tau}(\theta) \big)}{\tr\big( \bm{A}_{11}(\theta)\bm{A}_{11}(\theta) \big) + \tr\big( \bm{A}_{12}(\theta)\bm{A}_{21}(\theta) \big) +\tr\big( \bm{A}_{21}(\theta)\bm{A}_{12}(\theta) \big) +\tr\big( \bm{A}_{\tau}(\theta)\bm{A}_{\tau}(\theta) \big)},
	\end{equation}
	and
	\begin{equation} \label{lemma2_b}
	b = \frac{2i}{\tau} \frac{\frac{1}{2i\tau} \big[ \tr\big(\bm{A}_{\tau + i}(\theta)\bm{A}_{\tau + i}(\theta)\big) - \tr\big(\bm{A}_{\tau}(\theta)\bm{A}_{\tau}(\theta)\big) \big]}{\frac{1}{\tau^{2}} \tr\big(\bm{A}_{\tau}(\theta)\bm{A}_{\tau}(\theta)\big)}.
	\end{equation}
	
	As $i$ and $\delta$ are small compared to $\tau$ and $\theta$, the terms $\tr\big( \dot{\bm{A}}_{\tau}(\theta)\bm{A}_{\tau}(\theta) \big)$ and $\tr\big( \bm{A}_{\tau}(\theta)\bm{A}_{\tau}(\theta) \big)$ are relatively larger than the subdiagonal elements in \eqref{lemma2_a},
	and hence, $a$ can be further approximated as,
	\begin{equation} \label{lemma2_a1}
	a \approx  \frac{\tr\big( \dot{\bm{A}}_{\tau}(\theta)\bm{A}_{\tau}(\theta) \big)}{\tr\big( \bm{A}_{\tau}(\theta)\bm{A}_{\tau}(\theta) \big)}.
	\end{equation}
	
	Meanwhile, we approximate the term $\frac{1}{i\tau} \big[ \tr\big(\bm{A}_{\tau + i}(\theta)\bm{A}_{\tau + i}(\theta)\big) - \tr\big(\bm{A}_{\tau}(\theta)\bm{A}_{\tau}(\theta)\big) \big]$ in \eqref{lemma2_b} using $\frac{1}{\tau} \big[ \tr\big(\bm{A}_{\tau + 1}(\theta)\bm{A}_{\tau + 1}(\theta)\big) - \tr\big(\bm{A}_{\tau}(\theta)\bm{A}_{\tau}(\theta)\big) \big]$, and then we have,
	\begin{equation} \label{lemma2_b1}
	b \approx \frac{i}{\tau} \mu(\tau, \theta).
	\end{equation}
	
	The argument for the above approximation is as follows. First note that
	\begin{equation*}
	\begin{aligned}
	\bm{A}_{\tau + i}(\theta) &=  \bm{\Sigma}_{\tau+i}^{-1} \bm{V}_{\tau+i}(\theta) = (\bm{I}_{\tau+i} \otimes \bm{\Sigma})^{-1} (\bm{R}_{\tau+i}(\theta) \otimes \bm{\Lambda}) \\
	& = \bm{R}_{\tau+i}(\theta) \otimes (\bm{\Sigma}^{-1}\bm{\Lambda}).
	\end{aligned}
	\end{equation*}
	Then we have,
	\begin{equation*}
	\begin{aligned}
	\tr\big(\bm{A}_{\tau + i}(\theta)\bm{A}_{\tau + i}(\theta)\big) &= \tr \big( (\bm{R}_{\tau+i}(\theta) \otimes (\bm{\Sigma}^{-1}\bm{\Lambda}))(\bm{R}_{\tau+i}(\theta) \otimes (\bm{\Sigma}^{-1}\bm{\Lambda})) \big) \\
	&= \tr \big( (\bm{R}_{\tau+i}(\theta)\bm{R}_{\tau+i}(\theta)) \otimes (\bm{\Sigma}^{-1}\bm{\Lambda}\bm{\Sigma}^{-1}\bm{\Lambda}) \big) \\
	&= \tr \big( \bm{R}_{\tau+i}(\theta)\bm{R}_{\tau+i}(\theta) \big) \tr \big( \bm{\Sigma}^{-1}\bm{\Lambda}\bm{\Sigma}^{-1}\bm{\Lambda} \big) \\
	&= \tr \big( \bm{\Sigma}^{-1}\bm{\Lambda}\bm{\Sigma}^{-1}\bm{\Lambda} \big) \sum_{j} \sum_{k} [\bm{R}_{\tau+i}(\theta)]_{jk}^{2}  \\
	&=  \tr \big( \bm{\Sigma}^{-1}\bm{\Lambda}\bm{\Sigma}^{-1}\bm{\Lambda} \big) \Bigg(i \sum_{j} \sum_{k} [\bm{R}_{\tau+1}(\theta)]_{jk}^{2} + \sum_{|j - k|>\tau} [\bm{R}_{\tau+1}(\theta)]_{jk}^{2} \Bigg) \\
	& \approx  \tr \big( \bm{\Sigma}^{-1}\bm{\Lambda}\bm{\Sigma}^{-1}\bm{\Lambda} \big) \Bigg(i \sum_{j} \sum_{k} [\bm{R}_{\tau+1}(\theta)]_{jk}^{2}\Bigg).
	\end{aligned}
	\end{equation*}
	The last approximation is due to the fact that $(j, k)$th element of $\bm{R}_{\tau+1}(\theta)$ such that $|j - k|>\tau$ is small.

	Combining \eqref{lemma2_num}, \eqref{lemma2_den}, \eqref{lemma2_a1} \eqref{lemma2_b1} and the Taylor expansion $\frac{1}{\sqrt{1+x}} \approx 1 - \frac{1}{2}x + o(x)$, we obtain the approximation in \eqref{taylor1_main}.
	
\end{proof}

\end{document}